\RequirePackage{fix-cm}
\documentclass[smallextended]{svjour3}       
\smartqed  

\usepackage{amsfonts}
\usepackage{epstopdf}
\usepackage{latexsym}
\usepackage{amsmath}
\usepackage{setspace}
\usepackage{comment}
\usepackage{bm}
\usepackage{graphics}
\usepackage{graphicx}
\usepackage[caption=false]{subfig}
\usepackage{color}
\usepackage{algorithm}
\usepackage{algpseudocode}
\usepackage{algorithmicx}
\usepackage{amsopn}

\newcommand{\mini}{\mathop{\rm minimize}}
\newcommand{\maxi}{\mathop{\rm maximize}}

\begin{document}

\newtheorem{Ass}{Assumption}
\newtheorem{Al}{Algorithm}
\newtheorem{Def}{Definition}
\newtheorem{Rem}{Remark}
\newtheorem{pro}{Proposition}
\newtheorem{Cor}{Corollary}

\title{
A stabilized sequential quadratic semidefinite programming method for degenerate nonlinear semidefinite programs
\thanks{
The authors were in part supported by the Japan Society for the Promotion of Science KAKENHI for 20K19748, 20H04145, and 21K17709.
}
}

\titlerunning{A stabilized SQP method for nonlinear semidefinite programs}         

\author{Yuya Yamakawa \and Takayuki Okuno}


\institute{
Yuya Yamakawa \at Department of Applied Mathematics and Physics, Graduate School of Informatics, Kyoto University, Yoshidahommachi, Sakyo-ku, Kyoto-shi, Kyoto 606-8501, Japan,
\\
\email{yuya@i.kyoto-u.ac.jp}
\and
Takayuki Okuno \at The Center for Advanced Intelligence Project (AIP), RIKEN, Nihonbashi 1-chome Mitsui Building, 15th floor, 1-4-1 Nihonbashi, Chuo-ku, Tokyo 103-0027, Japan, 
\\
\email{takayuki.okuno.ks@riken.jp}
}

\date{Received: date / Accepted: date}

\maketitle

\begin{abstract}
In this paper, we propose a new sequential quadratic semidefinite programming (SQSDP) method for solving degenerate nonlinear semidefinite programs (NSDPs), in which we produce iteration points by solving a sequence of {\it stabilized} quadratic semidefinite programming (QSDP) subproblems, which we derive from the minimax problem associated with the NSDP. Unlike the existing SQSDP methods, the proposed one allows us to solve those QSDP subproblems inexactly, and each QSDP is feasible. One more remarkable point of the proposed method is that constraint qualifications (CQs) or boundedness of Lagrange multiplier sequences are not required in the global convergence analysis. Specifically, without assuming such conditions, we prove the global convergence to a point satisfying any of the following: the stationary conditions for the feasibility problem, the approximate-Karush-Kuhn-Tucker (AKKT) conditions, and the trace-AKKT conditions. Finally, we conduct some numerical experiments to examine the efficiency of the proposed method.
\end{abstract}

\section{Introduction}
In this paper, we consider the following nonlinear semidefinite program (NSDP):
\begin{equation}
\begin{array}{lll}
\displaystyle \mini_{x \in {\bf R}^{n}} \  & f(x)
\\ 
{\rm subject\ to}\  & g(x)=0, ~ X(x) \succeq O,
\end{array}
\label{NSDP}
\end{equation}
where $f \colon {\bf R}^{n} \rightarrow {\bf R}$, $g \colon {\bf R}^{n} \rightarrow {\bf R}^{m}$, and $X \colon {\bf R}^{n} \rightarrow {\bf S}^{d}$ are twice continuously differentiable functions, and ${\bf S}^{d}$ denotes the set of $d \times d$ real symmetric matrices. Let ${\bf S}_{++}^{d} \ ({\bf S}_{+}^{d})$ denote the set of $d \times d$ real symmetric positive (semi)definite matrices. For a matrix $M \in {\bf S}^{d}$, $M \succeq O$ and $M \succ O$ mean $M \in {\bf S}_{+}^{d}$ and $M \in {\bf S}_{++}^{d}$, respectively. Moreover, let $g_{1}, \ldots, g_{m}$ be the functions such that $g(x) = [g_{1}(x) \, \cdots \, g_{m}(x)]^{\top}$ for all $x \in {\bf R}^{n}$. In particular, when the functions $f$, $g$, and $X$ are affine, NSDP (\ref{NSDP}) is a linear semidefinite program (LSDP). If the range of $X$ is limited to the diagonal matrix space, then, it reduces to the standard nonlinear program (NLP).
\par
NSDPs have wide applications in control, finance, eigenvalue problems, structural optimization, and so forth \cite{St05,HYaHYa15,YaYaHa12}. Motivated by such practical importance, various algorithms have been developed for solving NSDPs so far, for example, sequential quadratic semidefinite programming (SQSDP) methods \cite{CoRa04,FaNoAp02,FrJaVo07,GoRa10,ZhCh16,zhao2018sqp,zhao2020line,ZhZh14}, interior point methods \cite{Ja00,LeMo02,okuno2020local,YaYa14,YaYa15,HYaHYa12,YaYaHa12}, augmented Lagrangian methods \cite{FuLo18,HuTeYa06,KoSt03,St05,SuSuZh08,WuLuDiCh13}, and others \cite{KaNaKaFu05,ZhAnSu13,yang2013homotopy}. However, in comparison with LSDPs and NLPs, studies on algorithms for NSDPs are still much fewer, and there is a lot of room for studying them in more depth. Our aim in this paper is to advance the SQSDP method further.
\par
Let us review existing works about the SQSDP method in more detail. The SQSDP method solves a sequence of quadratic SDP (QSDP) subproblems, which approximate the NSDP, so as to generate a sequence. It can be regarded as an extension of the sequential quadratic programming (SQP) method (e.g. see \cite{boggs1995sequential} for a survey) for the NLP. Strengths of the SQSDP method are that it admits both the global convergence and fast local convergence property, and moreover, it does not require a strictly feasible point as a starting point unlike interior point methods. We believe that the first SQSDP method was presented by Correa and Ram\'irez \cite{CoRa04} and the global convergence to a Karush-Kuhn-Tucker (KKT) point was established therein. The local convergence of the SQSDP was studied by Freund et al \cite{FrJaVo07}. Afterwards, some variants of SQSDP methods were proposed, e.g., a successive linearization method by Kanzow et al \cite{KaNaKaFu05}, a filter-type method by Zhu and Zhu \cite{ZhZh14}, and a penalty-free method equipped with the second-order correction step by Zhao and Chen \cite{ZhCh16}, who also proposed another SQSDP method that solves a feasible QSDP subproblem at every iteration in \cite{zhao2020line}.
\par
In this paper, we propose a new SQSDP method motivated from the stabilized SQP method for the NLP. The stabilized SQP method was initiated by Wright \cite{Wr98} to devise a superlinear convergent algorithm for solving degenerate NLPs. Since then, the stabilized SQP method have been developed by many researchers, e.g., \cite{GiKuRo17,gill2017stabilizedsup,GiRo13,Ha99,IzSoUs15,izmailov2017subspace}. Our SQSDP method, referred to as a stabilized SQSDP method, can be distinguished from the other SQSDP methods in that the following three points are achieved altogether:
\begin{description}
\item[(i)] {\it Consistent subproblem:} The proposed method organizes a consistent QSDP subproblem at each iteration, which always satisfies Slater's constraint qualification (CQ) and also has a unique optimum. 

\item[(ii)] {\it Inexact solution of QSDP:} We are allowed to truncate solving QSDP subproblems when certain criteria are met so as to ensure the global convergence.

\item[(iii)] {\it Stronger result on global convergence:} Without assuming boundedness of dual multiplier sequences or CQs at accumulation points, the global convergence is established successfully. In particular, this property is still valid even when NSDP~\eqref{NSDP} is infeasible.
\end{description}
We add some explanations regarding the above items. As for (i) and (ii), most the existing SQSDP methods are not well-defined in the sense that their subproblems are possibly infeasible. Moreover, they are impractical because exact optima of their subproblems are required to ensure the global convergence. Item~(iii) indicates another advantage of the proposed stabilized SQSDP method. All the existing SQSDP methods assume CQs or boundedness of a dual sequence, as far as we investigated. Those analyses may break down when applied to degenerate NSDPs, e.g., a convex NSDP such that Slater's CQ fails. In contrast, we establish the global convergence for the stabilized SQSDP method without such assumptions. More specifically, supposing the boundedness of a primal sequence generated by the stabilized SQSDP method, together with a certain controllable assumption, we prove that an accumulation point satisfies one of the following three conditions: a certain stationarity condition related to the feasibility, the AKKT conditions, and the TAKKT conditions. The AKKT and TAKKT conditions are the necessarily optimality conditions that were recently introduced by Andreani et al \cite{AnHaVi18}. The stabilized SQSDP method is expected to work effectively even when applied to degenerate or infeasible NSDPs.
\par
The proposed stabilized SQSDP method solves a QSDP subproblem that is derived from a regularized quadratic approximation to the minimax Lagrange problem for NSDP~\eqref{NSDP} to produce a search direction together with Lagrange multiplier estimates. If a point satisfying certain criteria is found, we may terminate solving the QSDP and proceed to the next step. A step size along the obtained search direction is determined with a backtracking line-search using the ordinary augmented Lagrangian. The next dual iterates and each parameter are determined according to the VOMF procedures that we will propose for NSDP~\eqref{NSDP}. The VOMF procedure determines an updating rule of parameters, which is rooted in the one presented by Gill and Robinson \cite{GiRo13} for solving NLPs. For the sake of clarity, we first present a prototype algorithm based on the VOMF procedure for solving the NSDP, and also give generic assumptions to obtain the global convergence. We next present the overall stabilized SQSDP method and prove the global convergence using the convergence result of the prototype algorithm. 
\par
The remaining part of this paper is organized as follows. In Section~\ref{sec_pre}, we introduce several notation and important concepts such as the KKT, AKKT, and TAKKT conditions. Section~\ref{sec:prototype_algorithm} provides a prototype algorithm for NSDP~(\ref{NSDP}) and its global convergence property. In Section~\ref{sec_SQSDP}, we propose a stabilized SQSDP method for finding an AKKT or a TAKKT point of NSDP~(\ref{NSDP}). Section~\ref{sec:convergence_analysis} proves the global convergence property of the stabilized SQSDP method. In Section~\ref{sec_numerical}, we report some numerical results associated with the stabilized SQSDP method. Finally, we make some concluding remarks in Section~\ref{sec_conclusion}.
\par
Throughout this paper, we use the following notation. For matrices $A$ and $B$ included in ${\bf R}^{p \times q}$, $\left\langle A, B \right\rangle$ represents the inner product of $A$ and $B$ defined by $\left\langle A, B \right\rangle := {\rm tr}(A^{\top}B)$, where ${\rm tr}(M)$ denotes the trace of a square matrix $M$, and the superscript $\top$ means the transposition of a vector or a matrix. Note that if $q = 1$, then $\left\langle \cdot, \cdot \right\rangle$ is the inner product of vectors in ${\bf R}^{p}$. The identity matrix and the all-ones vector are represented by $I$ and $e$, respectively, where these dimensions are defined by the context. For a vector $w \in {\bf R}^{p}$, $[w]_{i}$ indicates the $i$-th element of $w$, and $\Vert w \Vert$ denotes the Euclidean norm of $w$ defined by $\Vert w \Vert := \sqrt{\left\langle w, w \right\rangle}$. For a matrix $W \in {\bf R}^{p \times q}$, $[W]_{ij}$ is the $(i,j)$-th entry of $W$, and $\Vert W \Vert_{{\rm F}}$ means the Frobenius norm of $W$ defined by $\Vert W \Vert_{{\rm F}} := \sqrt{\left\langle W, W \right\rangle}$, and $\Vert W \Vert_{2}$ stands for the operator norm of $W$ defined by $\Vert W \Vert_{2} := \max \{ \Vert Wx \Vert \colon \Vert x \Vert = 1 \}$. For real numbers $r_{1}, \ldots, r_{d} \in {\bf R}$, we define
\begin{eqnarray*}
{\rm diag} \left[ r_{1}, \ldots, r_{d} \right]  := \left[
\begin{array}{ccc}
r_{1} & & O \\
&\ddots & \\
O & & r_{d}
\end{array}
\right].
\end{eqnarray*}
Let $U \in {\bf S}^{d}$ be a matrix with an orthogonal diagonalization $U = P D P^{\top}$, where $P$ is an orthogonal matrix and $D$ is a diagonal matrix. We denote by $\lambda_{1}^{P}(U), \ldots, \lambda_{d}^{P}(U)$ its eigenvalues satisfying $D = {\rm diag}[\lambda_{1}^{P}(U), \ldots, \lambda_{d}^{P}(U)]$. The minimum and the maximum eigenvalues of $U$ are expressed as $\lambda_{\min}(U)$ and $\lambda_{\max}(U)$, respectively. Furthermore, we denote by $[U]_{+}$ the projection of $U$ on ${\bf S}_{+}$, that is,
\begin{eqnarray*}
[U]_{+} := P {\rm diag} \left[ [\lambda_{1}^{P}(U)]_{+}, \ldots, [\lambda_{d}^{P}(U)]_{+} \right] P^{\top},
\end{eqnarray*}
where $[r]_{+} := \max \{ r, 0 \}$ for all $r \in {\bf R}$. Given open sets $P_{1}$, $P_{2}$, and $P_3$, let $\Phi$ be a mapping from $P_{1} \times P_{2}$ to $P_{3}$. We define the Fr\'echet derivative of $\Phi$ by $\nabla \Phi$. Moreover, we denote by $\nabla_{Z} \Phi$ the Fr\'echet derivative of $\Phi$ with respect to a variable $Z \in P_{1}$. For a closed convex set $S$, we write $\Pi_{S}$ for the metric projector over $S$. For a set $T$, the cardinality of $T$ is expressed as ${\rm card}(T)$. Finally, we will often use the following notation for the functions $g$ and $X$ in NSDP~\eqref{NSDP}:
\begin{itemize}
\item The matrix $\nabla g(x) \in {\bf R}^{n \times m}$ means $\nabla g(x) := [\nabla g_{1}(x) \, \cdots \, \nabla g_{m}(x)]$;
\item the matrix $A_{i}(x) \in {\bf S}^{d}$ indicates $A_{i}(x) := \frac{\partial}{\partial [x]_{i}} X(x)$ for $i=1, \ldots, n$;
\item the operator ${\cal A}(x) \colon {\bf R}^{n} \to {\bf S}^{d}$ is defined by ${\cal A}(x) u := [u]_{1} A_{1}(x) + \cdots + [u]_{n} A_{n}(x)$ for all $u \in {\bf R}^{n}$;
\item the adjoint operator of ${\cal A}(x)$ is represented by ${\cal A}^{\ast}(x) \colon {\bf S}^{d} \to {\bf R}^{n}$, that is, ${\cal A}^{\ast}(x) U = [\langle A_{1}(x), U \rangle \, \cdots \, \langle A_{n}(x), U \rangle]^{\top}$ for all $U \in {\bf S}^{d}$.
\end{itemize}

\section{Preliminaries} \label{sec_pre}
In this section, we define some notation and terminologies.

\subsection{The KKT conditions for NSDP~\eqref{NSDP}} \label{sectionKKT}
First, we introduce the Karush-Kuhn-Tucker (KKT) conditions for (\ref{NSDP}). Define the Lagrange function $L \colon {\bf R}^{n} \times {\bf R}^{m} \times {\bf S}^{d}_{+} \to {\bf R}$ as
\begin{eqnarray*}
L(v) := f(x) - \langle g(x), y \rangle - \left\langle X(x), Z \right\rangle,
\end{eqnarray*}
where $v := (x, y, Z)$. Note that $y \in {\bf R}^{m}$ and $Z \in {\bf S}^{d}_{+}$ are Lagrange multipliers for $g(x) = 0$ and $X(x) \succeq 0$, respectively. The gradient of $L$ at $v$ with respect to $x$ is given by 
\begin{eqnarray*}
\nabla _{x} L(v) = \nabla f(x) - \nabla g(x) y - {\cal A}^{\ast}(x) Z.
\end{eqnarray*}
The KKT conditions for (\ref{NSDP}) are represented in terms of the function $L$ as follows:

\begin{definition}
We say that $v = (x, y, Z) \in {\bf R}^{n} \times {\bf R}^{m} \times {\bf S}^{d}$ satisfies the KKT conditions if
\begin{eqnarray*} 
\displaystyle \nabla_{x} L(v) = 0, \quad g(x) = 0, \quad \langle X(x), Z \rangle = 0, \quad \displaystyle X(x) \succeq O, \quad Z \succeq O.
\end{eqnarray*}
In particular, we call the above point $x$ a KKT point and also call $(x, y, Z)$ a KKT triplet.
\end{definition}
The KKT conditions are the first-order optimality conditions for NSDP~\eqref{NSDP} under the presence of some CQ. We give two well-known CQs for NSDP~\eqref{NSDP}.
\begin{definition}
We say that a feasible point $x \in {\bf R}^{n}$ satisfies Robinson's CQ if
\begin{eqnarray*}
0 \in {\rm int} \left( \left[
\begin{array}{c}
g(x)
\\
X(x)
\end{array}
\right] + \left[
\begin{array}{c}
\nabla g(x)^{\top}
\\
{\cal A}(x)
\end{array}
\right] {\bf R}^{n} - \left[
\begin{array}{c}
\{ 0 \}
\\
{\bf S}^{d}_{+}
\end{array}
\right] \right),
\end{eqnarray*}
where ${\rm int}(S)$ denotes the topological interior of the set $S$. 
\end{definition}
\begin{definition}
We say that a feasible point $x \in {\bf R}^{n}$ satisfies the Mangasarian-Fromovitz constraint qualification (MFCQ) if
\begin{eqnarray}
& \{ \nabla g_{j}(x) \}_{j=1}^{m} ~ \mbox{{\rm are linearly independent,}} & \label{MFCQ1}
\\
& \exists d \in {\bf R}^{n} ~ \mbox{{\rm s.t.}} ~ \nabla g(x)^{\top} d = 0, ~ X(x) + {\cal A}(x) d \succ O. & \label{MFCQ2}
\end{eqnarray}
\end{definition}
A well-known fact is that, given a KKT point $x$, the corresponding Lagrange multipliers set $\{ (y,Z) \in {\bf R}^{m} \times {\bf S}^{d} \colon \mbox{$(x, y, Z)$ satisfies the KKT conditions} \}$ is nonempty and bounded under Robinson's CQ \cite{Ku76}. It follows from \cite[Chapter 3, prop. 2.1.12]{HiLe93} that Robinson's CQ is equivalent to the MFCQ. The global convergence properties of many existing methods are shown under the MFCQ or relevant assumptions.

\subsection{The AKKT and TAKKT conditions for NSDP~\eqref{NSDP}}
This section provides the definitions of the AKKT and TAKKT conditions for NSDP (\ref{NSDP}). These concepts have been introduced by Andreani, Haeser, and Viana \cite{AnHaVi18}. In the following, we first give the AKKT conditions:
\begin{definition} \label{AKKT_def}
We say that $x \in {\bf R}^{n}$ satisfies the AKKT conditions if $g(x)=0$, $X(x) \succeq O$, and there exist sequences $\{ x_{k} \} \subset {\bf R}^{n}, ~ \{ y_{k} \} \subset {\bf R}^{m}$, and $\{ Z_{k} \} \subset {\bf S}^{d}_{+}$ such that 
\begin{eqnarray*}
& \displaystyle \lim_{k \to \infty} x_{k} = x, \quad \displaystyle \lim_{k \to \infty} 
\nabla_{x} L(x_k,y_k,Z_k) = 0, &
\\
& \lambda_{j}^{U}(X(x)) > 0 \quad \Longrightarrow \quad \exists k_{j} \in {\bf N} \quad \mbox{{\rm s.t.}} \quad \lambda_{j}^{U_{k}}(Z_{k}) = 0 \quad \forall k \geq k_{j}, &
\end{eqnarray*}
where $j \in \{1, \ldots, d\}$, and $U$ and $U_{k} ~ (k \in {\bf N})$ are orthogonal matrices such that $U_{k} \to U ~ (k \to \infty)$,
\begin{eqnarray*}
& X(x) = U {\rm diag}[ \lambda_{1}^{U}(X(x)), \ldots, \lambda_{d}^{U}(X(x)) ] U^{\top}, &
\\
& Z_{k} = U_{k} {\rm diag}[ \lambda_{1}^{U_{k}}(Z_{k}), \ldots, \lambda_{d}^{U_{k}}(Z_{k}) ] U_{k}^{\top}. &
\end{eqnarray*}
\end{definition}
Moreover, we define the TAKKT conditions:
\begin{definition} \label{TAKKT_def}
We say that $x \in {\bf R}^{n}$ satisfies the TAKKT conditions if $g(x)=0$, $X(x) \succeq O$, and there exist sequences $\{ x_{k} \} \subset {\bf R}^{n}$, $\{ y_{k} \} \subset {\bf R}^{m}$, and $\{ Z_{k} \} \subset {\bf S}^{d}_{+}$ such that 
\begin{eqnarray*}
\displaystyle \lim_{k \to \infty} x_{k} = x, ~ \displaystyle \lim_{k \to \infty} \nabla_{x} L(x_k,y_k,Z_k) = 0, ~ \lim_{k \to \infty} \left\langle X(x_{k}), Z_{k} \right\rangle = 0.
\end{eqnarray*}
\end{definition}
In this paper, we call $x$ satisfying the AKKT conditions an AKKT point. Moreover, we call $\{ (x_{k}, y_{k}, Z_{k}) \}$ used for defining the AKKT point $x$ an AKKT sequence. As well, a TAKKT point and a TAKKT sequence are defined as for the TAKKT conditions.
\par
As was mentioned in the previous section, the KKT conditions are necessary optimality conditions under Robinson's CQ or the MFCQ. In contrast, the AKKT and TAKKT conditions always hold true as necessary optimality conditions in the absence of CQs, as stated in the next theorem.
\begin{theorem} {\rm \cite[Theorem 2, Theorem 5]{AnHaVi18}} \label{Th_nes}
Let $x$ be a local optimum of \eqref{NSDP}. Then, $x$ satisfies the AKKT and TAKKT conditions.
\end{theorem}
\par
The KKT conditions imply both the AKKT and TAKKT conditions. On the other hand, the AKKT conditions were shown not to imply the TAKKT conditions in \cite[Example 3]{AnHaVi18}. Moreover, the TAKKT conditions do not imply the AKKT conditions \cite[Example 3.1]{AnFuHaSaSe19}. Therefore, the AKKT and TAKKT conditions are mutually independent.
\par
The following theorem yields that if an AKKT or a TAKKT point satisfies the MFCQ, it is nothing but a KKT point.
\begin{theorem} {\rm \cite[Theorem 7, Theorem 8]{AnHaVi18}} \label{AKKTorTAKKTisKKT}
Let $x$ be a feasible point of NSDP \eqref{NSDP} satisfying the MFCQ. If $x$ satisfies the AKKT (TAKKT) conditions, then the AKKT (TAKKT) sequence corresponding to $x$ has a subsequence converging to a KKT point.
\end{theorem}
\par
Strictly speaking, the above AKKT and the TAKKT conditions together with Theorems~\ref{Th_nes} and \ref{AKKTorTAKKTisKKT} differ from the original ones presented in \cite{AnHaVi18}, in that the equality constraint $g(x) = 0$ is not handled there. However, by splitting $g(x) = 0$ into $g(x) \geq 0$ and $-g(x) \geq 0$, we can derive Theorems \ref{Th_nes} and \ref{AKKTorTAKKTisKKT} from \cite[Theorem 7, Theorem 8]{AnHaVi18}.

\section{A prototype algorithm} \label{sec:prototype_algorithm}
Before presenting the new SQSDP method, we provide its prototype algorithm and then study its convergence properties under general assumptions. This algorithm makes use of the following function $F \colon {\bf R}^{n} \to {\bf R}$ as a merit function for NSDP~\eqref{NSDP}:
\begin{eqnarray*} 
F(x; \sigma, y, Z) := f(x) + \frac{1}{2\sigma} \Vert \sigma y - g(x) \Vert^{2} + \frac{1}{2\sigma} \left\Vert [\sigma Z - X(x)]_{+} \right\Vert_{{\rm F}}^{2}.
\end{eqnarray*}
where $\sigma > 0$ is a penalty parameter. By \cite[Lemma 5]{AnHaVi18}, the function $F$ is continuously differentiable on ${\bf R}^{n}$, admitting the gradient
\begin{eqnarray}
\hspace{-3mm} \nabla F(x; \sigma, y, Z) = \nabla f(x) - \nabla g(x) \left\{ y - \frac{1}{\sigma} g(x) \right\} - {\cal A}^{\ast}(x) \left[ Z - \frac{1}{\sigma} X(x) \right]_{+}. \label{gradient_F}
\end{eqnarray}
Let ${\cal V} := {\bf R}^{n} \times {\bf R}^{m} \times {\bf S}^{d}$. For the later use, we also define $\widetilde{F} \colon {\cal V} \to {\bf R}$ by
\begin{eqnarray}
\widetilde{F}(x,y,Z; \sigma) := F(x; \sigma, y, Z) - \frac{\sigma}{2}(\Vert y \Vert^{2} + \Vert Z \Vert_{{\rm F}}^{2}), \label{augmented_LF}
\end{eqnarray}
which is the conventional augmented Lagrangian for NSDP~(\ref{NSDP}).
\subsection{Description of the prototype algorithm}
The prototype algorithm generates a sequence of primal iterates $\{ x_{k} \} \subset {\bf R}^{n}$ together with a sequence of Lagrange multiplier vectors and matrices $\{ (y_{k}, Z_{k}) \} \subset {\bf R}^{m} \times {\bf S}^{d}$, where $\{y_{k}\}$ and $\{ Z_{k} \}$ correspond to the constraints $g(x) = 0$ and $X(x) \succeq O$, respectively. Additionally, $\{ (\overline{y}_{k}, \overline{Z}_{k}) \} \subset {\bf R}^{m} \times {\bf S}^{d}$ is produced as a candidate for $\{ (y_{k}, Z_{k}) \}$.
\par
For each iteration $k \in {\bf N} \cup \{ 0 \}$, the algorithm attempts to solve
\begin{eqnarray} \label{mini_meritF0}
\begin{array}{ll}
\displaystyle \mini_{x \in {\bf R}^{n}} & F(x; \sigma_{k}, y_{k}, Z_{k}),
\end{array}
\end{eqnarray}
while tuning $(y_{k}, Z_{k})$ together with the penalty parameter $\sigma_{k} > 0$ according to a certain procedure which is clarified shortly (cf. Algorithm~\ref{ProcPara}). Although there still remain some parts not explained yet, let us show the overall figure of Algorithm~\ref{algorithm_SQSDP0} for the sake of understanding. The meaning of each step therein is as below. Hereafter, we denote the current iteration by $k \in {\bf N} \cup \{ 0 \}$ and define $v_{k} \in {\cal V}$ and $\overline{v}_{k} \in {\cal V}$ for every iteration $k$ as follows:
\begin{eqnarray*}
v_{k} := (x_{k}, y_{k}, Z_{k}) \quad \overline{v}_{k} := (x_{k}, \overline{y}_{k}, \overline{Z}_{k}).
\end{eqnarray*}
\par
In Step~1, $\overline{v}_{k+1} = (x_{k+1},\overline{y}_{k+1},\overline{Z}_{k+1})$ is output by approximately solving a problem relevant to (\ref{mini_meritF0}). We call this phase Mini-$F$-Phase. For the SQSDP method proposed in the next section, we will clarify how $\overline{v}_{k+1} = (x_{k+1}, \overline{y}_{k+1}, \overline{Z}_{k+1})$ is computed there. If $x_{k}$ is already a stationary point of  (\ref{mini_meritF0}), we set the Lagrange multiplier estimates $\overline{y}_{k+1}$ and $\overline{Z}_{k+1}$ in the spirit of the augmented Lagrangian method. Specifically, we put $\overline{y}_{k+1} := y_{k} - \frac{1}{\sigma_{k}} g(x_{k+1})$ and $\overline{Z}_{k+1} := [Z_{k} - \frac{1}{\sigma_{k}} X(x_{k+1})]_{+}$ with $x_{k+1} := x_{k}$. In Step~2, we update $(y_{k}, Z_{k})$ by performing the procedure, which is called VOMF-ITERATES and is described shortly, with $\overline{v}_{k+1} = (x_{k+1}, \overline{y}_{k+1}, \overline{Z}_{k+1})$ as input arguments. In Step~3, we decrease $\sigma_{k}$ as necessary to strengthen the penalty for the constraint violation in the function $F$. In what follows, we explain Steps~2 and 3 in more detail.

\subsection*{VOMF-ITERATES in Step~2}
VOMF-ITERATES is provided $\overline{v}_{k+1}$, $y_{k}$, $Z_{k}$, $\phi_{k}$, $\psi_{k}$, $\gamma_{k}$, and $\sigma_{k}$ as input arguments and generates a new Lagrange multiplier pair $(y_{k+1}, Z_{k+1})$ and a new parameter triplet $(\phi_{k+1}, \psi_{k+1}, \gamma_{k+1})$. It is formally described as in Algorithm~\ref{ProcPara}, wherein the following functions are utilized in order to measure the deviation of a given point $v = (x, y, Z) \in {\cal V}$ from the set of KKT points:
\begin{eqnarray*}
\Phi(v) := r_{V}(x) + \kappa r_{O}(v), \quad \Psi(v) := \kappa r_{V}(x) + r_{O}(v), \label{def_r1}
\end{eqnarray*}
where $\kappa \in (0,1)$ is a prefixed weight parameter and the functions $r_{V}$ and $r_{O}$ are defined by
\begin{eqnarray}
\hspace{5mm} r_{V}(x) := \Vert g(x) \Vert + [\lambda_{\max}(-X(x))]_{+}, \quad r_{O} (v) := \Vert \nabla_{x} L(v) \Vert + |\langle X(x), Z \rangle|. \label{def:rVrO}
\end{eqnarray}
Obviously, $v$ satisfies the KKT conditions if and only if $\Phi(v) = \Psi(v) = 0$. The procedure consists of the four steps, called the V-, O-, M-, and F-iterates, respectively. These names derive from those given in \cite{GiRo13} for nonlinear programming. In the V- and O- iterates, we check the KKT optimality of $\overline{v}_{k+1}$ with the value of $\Phi$ and $\Psi$. If $\Phi(\overline{v}_{k+1})$ and $\Psi(\overline{v}_{k+1})$ are not greater than $\frac{1}{2} \phi_{k}$ and $\frac{1}{2} \psi_{k}$, respectively, we regard $\overline{v}_{k+1}$ as a good approximation to KKT point, and set $y_{k+1} := \overline{y}_{k+1}$ and $Z_{k+1} := \overline{Z}_{k+1}$. We then decrease $\phi_{k}$ or $\psi_{k}$ to gain a more refined point in subsequent iterations. Otherwise, we proceed to the M-iterate which examines whether the following inequality holds or not:
\begin{eqnarray}
\Vert \nabla F(x_{k+1}; \sigma_{k}, y_{k}, Z_{k}) \Vert \leq \gamma_{k}. \label{AKKTstep}
\end{eqnarray}
If it holds true, we regard $x_{k+1}$ as a good approximate solution of
\begin{eqnarray} \label{mini_meritF}
\begin{array}{ll}
\displaystyle \mini_{x \in {\bf R}^{n}} & F(x; \sigma_{k}, y_{k}, Z_{k}).
\end{array}
\end{eqnarray}
In fact, the above problem~\eqref{mini_meritF} is identical to the subproblem of the augmented Lagrangian method using the augmented Lagrangian $\widetilde{F}$ defined by (\ref{augmented_LF}). The following updating formulae in the M-iterate are motivated by this fact:
\begin{eqnarray*}
\textstyle y_{k+1} := \Pi_{C} ( y_{k} - \frac{1}{\sigma_{k}}g(x_{k+1}) ), \quad Z_{k+1} := \Pi_{D} ( [ Z_{k} - \frac{1}{\sigma_{k}}X(x_{k+1}) ]_{+} ),
\end{eqnarray*}
where $C \subset {\bf R}^{m}$ and $D \subset {\bf S}^{d}$ are closed bounded sets defined in Algorithm~\ref{ProcPara}, and further decrease $\gamma_{k}$ so as to find a solution satisfying $\nabla F(x_{k}; \sigma_{k}, y_{k}, Z_{k}) = 0$ more accurately through subsequent iterations.
\par
When the if-statement related to the M-iterate is false, then we proceed to the F-iterate, therein just incrementing $k$ by one without any other updates. This step is aimed at gaining a refined point such that one of the criteria for the V-, O-, and M-iterates is satisfied.
\par
Roughly speaking, the V- and O- iterates play a part of generating a TAKKT sequence, while the M-iterate is for an AKKT sequence.

\subsection*{Update of the parameter $\sigma_{k}$ in Step~3}
If the condition (\ref{AKKTstep}) is fulfilled, we decrease $\sigma_{k}$ to strengthen the penalty for the constraint violation in the function $F$. Specifically, we update $\sigma_{k}$ as follows:
\begin{eqnarray} \label{update_sigma} 
\sigma_{k+1} := \left\{
\begin{array}{ll} 
\min \{ \frac{1}{2} \sigma_{k}, r(v_{k+1})^{\frac{3}{2}} \} & \mbox{if (\ref{AKKTstep}) is satisfied},
\\
\sigma_{k} & \mbox{otherwise}.
\end{array}
\right.
\end{eqnarray}
Although the term $r(v_{k+1})^{\frac{3}{2}}$ is employed in the above aiming at fast local convergence, other update formulae can be adopted as long as $\{ \sigma_{k} \}$ is a monotonically decreasing sequence of positive numbers.

\begin{algorithm}[tbh]
\caption{Prototype algorithm}\label{algorithm_SQSDP0}
\begin{algorithmic}[1]
\Require 
Choose $v_{0} := (x_{0}, y_{0}, Z_{0})$ such that $Z_{0} \succeq O$. Set
\begin{eqnarray*}
k := 0, ~ \overline{y}_{0} := y_{0}, ~ \overline{Z}_{0} := Z_{0}, ~ \phi_{0} > 0, ~ \psi_{0} > 0, ~ \gamma_{0} > 0, ~ \sigma_{0} > 0.
\end{eqnarray*}
Set parameters for Mini-$F$-Phase and {\rm VOMF-ITERATES} (see Lines~\ref{line4} and \ref{line5}).  
\Repeat
\If{$\Vert \nabla F(x_{k}; \sigma_{k}, y_{k}, Z_{k}) \Vert = 0$}
\Comment{Step~1~(Mini-$F$-Phase)} \label{line4}
\State{Set
\begin{eqnarray*}
\qquad \textstyle x_{k+1} := x_{k}, ~ \overline{y}_{k+1} := y_{k} - \frac{1}{\sigma_{k}} g(x_{k+1}), ~ \overline{Z}_{k+1} := [ Z_{k} - \frac{1}{\sigma_{k}} X(x_{k+1}) ]_{+}.
\end{eqnarray*}
}
\Else
\State{
Compute $\overline{v}_{k+1} = (x_{k+1},\overline{y}_{k+1},\overline{Z}_{k+1})$ by approximately solving \eqref{mini_meritF} or its relevant problem (cf. Lines~\ref{st:chooH}--\ref{st:xk} in Algorithm~\ref{algorithm_SQSDP}). 
}
\EndIf
\State{Set \Comment{Step~2~(VOMF-ITERATES)} \label{line5}
\begin{eqnarray*}
&& \hspace{-3mm} (y_{k+1},Z_{k+1},\phi_{k+1},\psi_{k+1},\gamma_{k+1})
\\
&& \hspace{19mm} := \mbox{{\rm VOMF-ITERATES}}(\overline{v}_{k+1}, y_{k}, Z_{k}, \phi_{k}, \psi_{k}, \gamma_{k}, \sigma_{k}).
\end{eqnarray*}
}
\State{
Update $\sigma_{k}$ by \Comment{Step~3}
\begin{eqnarray*} 
&& \sigma_{k+1} := \left\{
\begin{array}{ll}
\min \{ \frac{1}{2} \sigma_{k}, r(v_{k+1})^{\frac{3}{2}} \} & {\rm if ~} \Vert \nabla F(x_{k+1}; \sigma_{k}, y_{k}, Z_{k}) \Vert \leq \gamma_{k},
\\
\sigma_{k} & {\rm otherwise}.
\end{array}
\right.
\end{eqnarray*}}
\State{Set $k := k+1$.}
\Comment{Step~4}

\Until $v_{k} := (x_{k}, y_{k}, Z_{k})$ meets a suitable criterion
\end{algorithmic}
\end{algorithm}

\begin{algorithm}[tbh]
\caption{Procedure for updating $(y_{k}, Z_{k}, \phi_{k}, \psi_{k}, \gamma_{k})$} \label{ProcPara}
\begin{algorithmic}[1]
\Procedure{{\rm VOMF-ITERATES}}{$\overline{v}_{k+1}, y_{k}, Z_{k}, \phi_{k}, \psi_{k}, \gamma_{k}, \sigma_k$}

\If{$\Phi(\overline{v}_{k+1}) \leq \frac{1}{2} \phi_k$}
\State{
Set \Comment{V-iterate}
\begin{eqnarray*}
\textstyle \qquad y_{k+1} := \overline{y}_{k+1}, ~ Z_{k+1} := \overline{Z}_{k+1}, ~ \phi_{k+1} := \frac{1}{2} \phi_{k}, ~ \psi_{k+1} := \psi_{k}, ~ \gamma_{k+1} := \gamma_{k}.
\end{eqnarray*}
}

\ElsIf{$\Psi(\overline{v}_{k+1}) \leq \frac{1}{2} \psi_k$}
\State{
Set \Comment{O-iterate}
\begin{eqnarray*}
\textstyle \qquad y_{k+1} := \overline{y}_{k+1}, ~ Z_{k+1} := \overline{Z}_{k+1}, ~ \phi_{k+1} := \phi_{k}, ~ \psi_{k+1} := \frac{1}{2} \psi_{k}, ~ \gamma_{k+1} := \gamma_{k},
\end{eqnarray*}
}

\ElsIf{$\Vert \nabla F(x_{k+1}; \sigma_{k}, y_{k}, Z_{k}) \Vert \leq \gamma_{k}$}
\State{
Set \Comment{M-iterate}
\begin{eqnarray*}
& \textstyle y_{k+1} := \Pi_{C}(y_{k} - \frac{1}{\sigma_{k}}g(x_{k+1})), ~ Z_{k+1} := \Pi_{D}([Z_{k} - \frac{1}{\sigma_{k}} X(x_{k+1})]_{+}), &
\\
& \textstyle \phi_{k+1} := \phi_{k},\ \psi_{k+1} := \psi_{k}, ~ \gamma_{k+1} := \frac{1}{2}\gamma_{k}, &
\end{eqnarray*}
where $C := \{y \in {\bf R}^{m} \colon -y_{\max}e \leq y\leq y_{\max} e \}$, $D := \{ Z\in {\bf S}^{d} \colon O \preceq Z \preceq z_{\max} I \}$.
}

\Else
\State{
Set \Comment{F-iterate}
\begin{eqnarray*}
y_{k+1} := y_{k}, ~ Z_{k+1} := Z_{k}, ~ \phi_{k+1} := \phi_{k}, ~ \psi_{k+1} := \psi_{k}, ~ \gamma_{k+1} := \gamma_{k}.
\end{eqnarray*}
}
\EndIf 
\State{{\bf return} $(y_{k+1},Z_{k+1},\phi_{k+1},\psi_{k+1},\gamma_{k+1})$} 
\EndProcedure
\end{algorithmic}
\end{algorithm}

\subsection{Convergence analysis of Algorithm~\ref{algorithm_SQSDP0}} \label{sec:prototype_global_convergence}
We prove the global convergence property of Algorithm~\ref{algorithm_SQSDP0}. In the subsequent arguments, we use three sets ${\cal I}$, ${\cal J}$, and ${\cal K}$ defined by
\begin{eqnarray} \label{def:sets_IJK}
\begin{array}{rcl}
\hspace{-5mm} {\cal I} \hspace{0mm} &:=& \hspace{0mm} \{ k \colon y_{k}, \, Z_{k}, \, \phi_{k}, \, \psi_{k}, \, \mbox{and} ~ \gamma_{k} ~ \mbox{are updated by the V- or O-iterate} \}, 
\\
\hspace{-5mm}{\cal J} \hspace{0mm} &:=& \hspace{0mm} \{ k \colon y_{k}, \, Z_{k}, \, \phi_{k}, \, \psi_{k}, \, \mbox{and} ~ \gamma_{k} ~ \mbox{are updated by the M-iterate} \},
\\
\hspace{-5mm}{\cal K} \hspace{0mm} &:=& \hspace{0mm} \{ k \colon y_{k}, \, Z_{k}, \, \phi_{k}, \, \psi_{k}, \, \mbox{and} ~ \gamma_{k} ~ \mbox{are updated by the F-iterate} \}.
\end{array}
\end{eqnarray}
Notice that ${\cal I}$, ${\cal J}$, and ${\cal K}$ are mutually disjoint. Moreover, we make the following two sets of assumptions:
\begin{Ass} \label{assumption_global0}
\
\begin{description}
\item[\rm (A1)] The functions $f$, $g$, and $X$ are twice continuously differentiable;
\item[\rm (A2)] there exists a compact set $\Gamma$ such that any sequence $\{ x_{k} \}$ generated by Algorithm~{\rm \ref{algorithm_SQSDP0}} is contained in $\Gamma$;
\end{description}
\end{Ass}
\begin{Ass} \label{assumption_IJK}
There never occurs a situation with ${\rm card}({\cal I}) < \infty$, ${\rm card}({\cal J}) < \infty$, and ${\rm card}({\cal K}) = \infty$.
\end{Ass}

Assumption~\ref{assumption_global0} is standard and can be found in many literatures. It is also assumed in the convergence analysis for a stabilized SQSDP method, which is presented in the next section, whereas Assumption~\ref{assumption_IJK} will be verified there. In addition to the above, we also implicitly assume that Algorithm~\ref{algorithm_SQSDP0} generates infinitely many iteration points.
\par
In the following, we give properties associated with $\{ \phi_{k} \}$, $\{ \psi_{k} \}$, $\{ \gamma_{k} \}$, $\{ \sigma_{k} \}$, and $\{ Z_{k} \}$. Its proof is given in Appendix~\ref{app_lemindex}.

\begin{lemma} \label{index_lemma}
Suppose that Assumption~{\rm \ref{assumption_global0}} holds. Then, we have the following properties:
\begin{description}
\item[{\rm (i)}] If ${\rm card}({\cal I}) = \infty$, then $\phi_{k} \to 0$ or $\psi_{k} \to 0$ as $k \to \infty$;
\item[{\rm (ii)}] if ${\rm card}({\cal I}) < \infty$, then $\{Z_{k} \}$ is bounded;
\item[{\rm (iii)}] if ${\rm card}({\cal I}) < \infty$ and ${\rm card}({\cal J}) = \infty$, then $\sigma_{k} \to 0$ and $\gamma_{k} \to 0$ as $k \to \infty$.
\end{description}
\end{lemma}

With the help of the above lemma, we derive the following convergence theorem of Algorithm~\ref{algorithm_SQSDP0}.

\begin{theorem} \label{global_convergence}
Suppose that Assumptions~{\rm \ref{assumption_global0}} and {\rm \ref{assumption_IJK}} hold. Any accumulation point of $\{ x_{k} \}$, say $x^{\ast}$, satisfies at least one of the following statements:
\begin{description}
\item[{\rm (i)}] $x^{\ast}$ is a TAKKT point of \eqref{NSDP};
\item[{\rm (ii)}] $x^{\ast}$ is an AKKT point of \eqref{NSDP};
\item[{\rm (iii)}] $x^{\ast}$ is an infeasible point of \eqref{NSDP}, but a stationary point of the following optimization problem:
\begin{eqnarray*}
\begin{array}{ll}
\displaystyle \mini_{x \in {\bf R}^{n}} & h(x) := \displaystyle \frac{1}{2} \Vert g(x) \Vert^{2} + \frac{1}{2} \Vert [-X(x)]_{+} \Vert_{{\rm F}}^{2},
\end{array}
\end{eqnarray*}
that is to say, $\nabla h(x^{\ast}) = 0$.
\end{description}
\end{theorem}

\noindent
{\it Proof.}
We consider the two cases: (a) ${\rm card}({\cal I}) = \infty$; (b) ${\rm card}({\cal I}) < \infty$. 
\\
{\bf Case (a):} Let ${\cal P} := \{ k \in {\bf N} \colon k-1 \in {\cal I} \}$. Note that ${\rm card}({\cal P}) = \infty$ according to ${\rm card}({\cal I}) = \infty$. Recall that $v_{k} = (x_{k}, y_{k}, Z_{k})$ and $\overline{v}_{k} = (x_{k}, \overline{y}_{k}, \overline{Z}_{k})$. Assumption~\ref{assumption_global0} (A2) implies that $\{ x_{k} \}_{{\cal P}}$ has at least one accumulation point, say $x^{\ast}$. Then, there exists ${\cal L} \subset {\cal P}$ such that $x_{k} \to x^{\ast}$ as ${\cal L} \ni k \to \infty$. It follows from Lemma \ref{index_lemma} (i) that $\phi_{k} \to 0$ or $\psi_{k} \to 0$ as ${\cal L} \ni k \to \infty$. Since $\Phi(v_{k}) = \Phi(\overline{v}_{k}) \leq \frac{1}{2}\phi_{k-1} = \phi_{k}$ or $\Psi(v_{k}) = \Psi(\overline{v}_{k}) \leq \frac{1}{2}\psi_{k-1} = \psi_{k}$ for $k \in {\cal L} \subset {\cal P}$, it is clear that $\Phi(v_{k}) \to 0$ or $\Psi(v_{k}) \to 0$ as ${\cal L} \ni k \to \infty$, that is,
\begin{eqnarray*}
& \displaystyle \lim_{{\cal L} \ni k \to \infty} \nabla_{x} L(x_{k}, y_{k}, Z_{k}) = 0, \, \lim_{{\cal L} \ni k \to \infty} \langle Z_{k}, X(x_{k}) \rangle = 0, & 
\\ 
& \hspace{-0.5mm} g(x^{\ast}) = \displaystyle \lim_{{\cal L} \ni k \to \infty} g(x_{k}) = 0, \, [\lambda_{\max}(-X(x^{\ast}))]_{+} = \lim_{{\cal L} \ni k \to \infty} [\lambda_{\max}(-X(x_{k}))]_{+} = 0. & 
\end{eqnarray*} 
These results and $\{ Z_{k} \} \subset {\bf S}^{d}_{+}$ derive that $\{ (x_{k}, y_{k}, Z_{k}) \}_{{\cal L}}$ is a TAKKT sequence corresponding to $x^{\ast}$. Hence, in this case, we have situation (i).
\\
{\bf Case (b):} We show that, in this case, there occurs (ii) or (iii). Let ${\cal Q} := \{ k \in {\bf N} \colon k-1 \in {\cal J} \}$, $\widetilde{y}_{k} := y_{k-1} - \frac{1}{\sigma_{k-1}}g(x_{k})$, and $\widetilde{Z}_{k} := [Z_{k-1} - \frac{1}{\sigma_{k-1}}X(x_{k}) ]_{+}$. Assumption~\ref{assumption_IJK} implies that ${\rm card}({\cal J}) = \infty$ must hold. Notice that ${\rm card}({\cal Q}) = \infty$ by ${\rm card}({\cal J}) = \infty$. Moreover, Assumption~\ref{assumption_global0}~(A2) ensures that $\{ x_{k} \}_{{\cal Q}}$ has at least one accumulation point, say $x^{\ast}$ again. Then, there exists ${\cal M} \subset {\cal Q}$ such that $x_{k} \to x^{\ast}$ as ${\cal M} \ni k \to \infty$. We have from Lemma \ref{index_lemma}~(iii) that $\sigma_{k-1} \to 0$ and $\gamma_{k-1} \to 0$ as ${\cal M} \ni k \to \infty$. Note that $\nabla_{x} L(x_{k}, \widetilde{y}_{k}, \widetilde{Z}_{k}) = \nabla f(x_{k}) - \nabla g(x_{k}) \widetilde{y}_{k} - {\cal A}^{\ast}(x_{k}) \widetilde{Z}_{k} = \nabla F(x_{k}; \sigma_{k-1}, y_{k-1}, Z_{k-1})$ and $\Vert \nabla F(x_{k}; \sigma_{k-1}, y_{k-1}, Z_{k-1}) \Vert \leq \gamma_{k-1}$ for $k \in {\cal M} \subset {\cal Q}$. These facts yield that
\begin{eqnarray}
& \displaystyle \lim_{{\cal M} \ni k \to \infty} \nabla_{x} L(x_{k}, \widetilde{y}_{k}, \widetilde{Z}_{k}) = 0, & \label{AKKT_1}
\\
& \displaystyle \nabla h(x^{\ast}) = \lim_{{\cal M} \ni k \to \infty} \sigma_{k-1} \nabla F(x_{k}; \sigma_{k-1}, y_{k-1}, Z_{k-1}) = 0. & \label{stationary_cond}
\end{eqnarray}
Furthermore, since the boundedness of $\{ Z_{k} \}$ is ensured by Lemma~\ref{index_lemma}~(ii), we obtain
\begin{eqnarray}
\displaystyle \lim_{{\cal M} \ni k \to \infty} \left\{ \sigma_{k-1} Z_{k-1} - X(x_{k}) \right\} = - X(x^{\ast}). \label{AKKT_XZ}
\end{eqnarray}
Denote $X^{\ast} := X(x^{\ast}), ~ X_{k} := X(x_{k})$, and $W_{k} := \sigma_{k-1} Z_{k-1} - X(x_{k})$. Take a diagonal decomposition of $X^{\ast}$, i.e., $X^{\ast} = V^{\ast} {\rm diag}[ \lambda_{1}^{V^{\ast}}(X^{\ast}), \ldots, \lambda_{d}^{V^{\ast}}(X^{\ast}) ] (V^{\ast})^{\top}$, where $V^{\ast}$ is an orthogonal matrix. By (\ref{AKKT_XZ}), there exists ${\cal N} \subset {\cal M}$ such that
\begin{eqnarray}
\lambda_{1}^{V_{k}}(W_{k}) \to - \lambda_{1}^{V^{\ast}}(X^{\ast}), \ldots, \lambda_{d}^{V_{k}}(W_{k}) \to - \lambda_{d}^{V^{\ast}}(X^{\ast}) \quad ({\cal N} \ni k \to \infty), \label{lambdato0}
\end{eqnarray}
where $V_{k} ~ (k \in {\cal L})$ are orthogonal matrices such that $V_{k} \to V^{\ast} ~ ({\cal N} \ni k \to \infty)$ and $W_{k} = V_{k}  {\rm diag}[ \lambda_{1}^{V_{k}}(W_{k}), \ldots, \lambda_{d}^{V_{k}}(W_{k}) ] V_{k}^{\top}$. Let $j \in \{ 1, \ldots, d \}$ be an arbitrary integer. Now, assume that $\lambda_{j}^{V^{\ast}}(X^{\ast}) > 0$. Combining this assumption and (\ref{lambdato0}) gives $\lim_{{\cal N} \ni k \to \infty} \lambda_{j}^{V_{k}}(W_{k}) = - \lambda_{j}^{V^{\ast}}(X^{\ast}) < 0$, and hence there exists $k_{j} \in {\bf N}$ such that $\lambda_{j}^{V_{k}}(W_{k}) < 0$ for all $k \in \{ k \in {\cal N} \colon k \geq k_{j} \}$. This fact and $\widetilde{Z}_{k} = [Z_{k-1} - \frac{1}{\sigma_{k-1}}X_{k}]_{+} = [\frac{1}{\sigma_{k-1}}W_{k}]_{+}$ yield $\lambda_{j}^{V_{k}}(\widetilde{Z}_{k}) = \frac{1}{\sigma_{k-1}} [\lambda_{j}^{V_{k}} (W_{k})]_{+} = 0$ for all $k \in \{ k  \in {\cal N} \colon k \geq k_{j} \}$. As a result, we obtain
\begin{eqnarray}
\hspace{-5mm} \lambda_{j}^{V^{\ast}} (X^{\ast}) > 0 ~~ \Longrightarrow ~~ \exists k_{j} \in {\bf N} ~~ {\rm s.t.} ~~ \lambda_{j}^{V_{k}} (\widetilde{Z}_{k}) = 0 ~~ \forall k \in \{ k  \in {\cal N} \colon k \geq k_{j} \}. \label{AKKT_2}
\end{eqnarray}
If $x^{\ast}$ is feasible to NSDP (\ref{NSDP}), it follows from (\ref{AKKT_1}) and (\ref{AKKT_2}) that $\{ (x_{k}, \widetilde{y}_{k}, \widetilde{Z}_{k}) \}_{{\cal N}}$ is an AKKT sequence corresponding to $x^{\ast}$, that is, situation (ii) holds. Even if $x^{\ast}$ is not feasible, situation (iii) holds from (\ref{stationary_cond}).
\par
Therefore, we conclude that if Case (a) holds, then situation (i) occurs, otherwise, namely, if Case (b) holds, situation (ii) or (iii) is satisfied.
$ \hfill \Box $
\bigskip

Finally, Theorem\,\ref{global_convergence} entails the following one about global convergence to a KKT point under the presence of the conditions related to the MFCQ.

\begin{theorem}\label{global_convergence2}
Suppose that Assumptions~{\rm \ref{assumption_global0}} and {\rm \ref{assumption_IJK}} hold. If any accumulation point of $\{ x_{k} \}$, say $x^{\ast}$, satisfies \eqref{MFCQ1} and \eqref{MFCQ2} in the definition of the MFCQ, then $x^{\ast}$ is nothing but a KKT point.
\end{theorem}

\noindent
{\it Proof.}
Note that AKKT and TAKKT points are KKT points if they satisfy conditions (\ref{MFCQ1}) and (\ref{MFCQ2}). Hence, to prove the desired claim, it is sufficient to show that $x^{\ast}$ is an AKKT or a TAKKT point. We prove this assertion by contradiction, that is, suppose that $x^{\ast}$ is neither an AKKT nor a TAKKT point. By Theorem \ref{global_convergence}, we see that $x^{\ast}$ is infeasible and $\nabla h(x^{\ast}) = 0$, yielding
\begin{eqnarray}
& \frac{1}{2} \Vert g(x^{\ast}) \Vert^{2} + \frac{1}{2} \Vert [-X(x^{\ast})]_{+} \Vert_{{\rm F}}^{2} = h(x^{\ast}) \not = 0, & \label{feasible_cond}
\\
& \nabla g(x^{\ast})g(x^{\ast}) - \mathcal{A}^{\ast}(x^{\ast})[-X(x^{\ast})]_{+} = \nabla h(x^{\ast}) = 0. & \label{hx0}
\end{eqnarray}
Since conditions (\ref{MFCQ1}) and (\ref{MFCQ2}) are satisfied at $x^{\ast}$, the matrix $\nabla g(x^{\ast})$ is of full rank, and there exists $d \in {\bf R}^{n}$ such that $\nabla g(x^{\ast})^{\top} d = 0$ and $X(x^{\ast}) + {\cal A}(x^{\ast}) d \succ O$. It then follows from $[-X(x^{\ast})]_{+} \succeq O$ that $0 \leq \langle [-X(x^{\ast})]_{+}, X(x^{\ast}) + {\cal  A}(x^{\ast})d \rangle$, which implies
\begin{eqnarray}
&& 0 \geq 
\langle [-X(x^{\ast})]_{+}, -X(x^{\ast}) \rangle - \langle {\cal  A}^{\ast}(x^{\ast})[-X(x)]_{+}, d \rangle \nonumber
\\
&& \phantom{0} = \Vert [-X(x^{\ast})]_{+} \Vert_{{\rm F}}^{2} - \langle \nabla g(x^{\ast}) g(x^{\ast}), d \rangle \nonumber
\\
&& \phantom{0} = \Vert [-X(x^{\ast})]_{+} \Vert_{{\rm F}}^{2} - \langle g(x^{\ast}), \nabla g(x^{\ast})^{\top} d \rangle, \label{Xgineq}
\end{eqnarray}
where the first equality is obtained by (\ref{hx0}). Combining (\ref{Xgineq}) with $\nabla g(x^{\ast})^{\top} d = 0$ yields that $\Vert [-X(x^{\ast})]_{+} \Vert_{{\rm F}} = 0$. Then, we have from (\ref{hx0}) that $\nabla g(x^{\ast}) g(x^{\ast}) = 0$. Since $\nabla g(x^{\ast})$ is of full rank, we get $g(x^{\ast}) = 0$. Therefore, $h(x^{\ast}) = 0$. However, this contradicts (\ref{feasible_cond}). The proof is complete. $ \hfill \Box $

\section{An inexact and stabilized SQSDP method} \label{sec_SQSDP}
In this section, we propose a new SQSDP method for NSDP (\ref{NSDP}) by integrating a sequential quadratic optimization technique into the prototype algorithm, i.e., Algorithm~\ref{algorithm_SQSDP0}. Specifically, we clarify how Mini-$F$-Phase of the algorithm is performed. At every Mini-F-Phase, a convex quadratic SDP (QSDP) is solved to compute a search direction $p$, and then a step size $\alpha > 0$ is determined along $p$ by means of a line-search technique equipped with the function $F$ as a merit function. We then update $x$ as $x + \alpha p$.
\par
Differently from the existing SQSDP methods, the proposed one may truncate a process of solving the convex QSDP subproblem if certain conditions are satisfied. Another remarkable point is that the QSDP is {\it strictly} feasible (Slater's CQ is fulfilled) and thus the KKT conditions always hold at its optimum.
\par
Hereafter, a current iterate is often represented as $(x,y,Z)$ simply.

\subsection{A stabilized quadratic SDP subproblem}
This section specifies how Mini-$F$-Phase is conducted. We first derive the quadratic SDP subproblem to be solved at $(x, y, Z)$ so as to minimize the merit function $F$. As shown in the subsequent argument, an optimal solution of this problem is indeed the descent direction of the merit function $F$ under some reasonable assumptions.
\par
We shall start with the Lagrange minimax problem regarding NSDP~(\ref{NSDP}):
\begin{eqnarray*}
\mini_{x\in {\bf R}^{n}} \maxi_{(y,Z)\in {\bf R}^{m} \times {\bf S}^{d}_{+}} L(x, y, Z).
\end{eqnarray*}
For solving this problem, similarly to Wright \cite{Wr98}, we iteratively solve the following problem which approximates the above problem at $(x, y, Z)$:
\begin{eqnarray} \label{minimax}
\begin{array}{ll}
\displaystyle \hspace{-5mm} \mini_{\xi \in {\bf R}^{n}} \maxi_{(\zeta, \Sigma)\in {\bf R}^{m} \times {\bf S}^{d}_{+}} & \displaystyle \langle\nabla f(x),\xi \rangle +\frac{1}{2} \langle H \xi, \xi \rangle - \langle \zeta, g(x) + \nabla g(x)^{\top}\xi \rangle
\\
& \displaystyle \hspace{0mm} - \langle \Sigma, X(x) + {\cal A}(x) \xi \rangle - \frac{\sigma}{2} \Vert \zeta - y \Vert^{2} - \frac{\sigma}{2} \Vert \Sigma - Z \Vert_{{\rm F}}^{2},
\end{array}
\end{eqnarray}
where $\sigma>0$ is a parameter shared with the function $F$ and $H \in {\bf R}^{n \times n}$ can be set to any matrix as long as Assumption~\ref{assumption_global4}, which will appear later, is satisfied for the global convergence. Although an arbitrary positive definite matrix is theoretically accepted as $H$, the Hessian $\nabla_{xx}^{2} L(x, y, Z)$ or its approximation is preferable for the sake of rapid convergence. Notice that the maximization part in (\ref{minimax}) is a strongly concave problem by virtue of the last two quadratic terms, thus it always has a unique maximum $(\zeta,\Sigma)$ for a prefixed $\xi$. Notice also that an optimum $\xi$ of (\ref{minimax}) together with the corresponding maximizer $(\zeta, \Sigma)$ solves the following system:
\begin{eqnarray}
& \nabla f(x) + H \xi - \nabla g(x) \zeta - {\cal A}^{\ast}(x) \Sigma = 0, & \label{KKT_minimax1}
\\
& g(x)+\nabla g(x)^{\top} \xi + \sigma ( \zeta - y ) = 0, & \label{KKT_minimax2}
\\
& {\bf S}^{d}_{+} \ni X(x) + {\cal A}(x) \xi + \sigma ( \Sigma - Z ) \perp \Sigma \in {\bf S}^{d}_{+}, & \label{KKT_minimax3}
\end{eqnarray}
where $\perp$ stands for the perpendicularity in ${\bf S}^{d}$. If $H$ is positive definite, then (\ref{KKT_minimax1})--(\ref{KKT_minimax3}) are necessary and sufficient conditions under which $(\xi, \zeta, \Sigma)$ is a unique optimal solution of (\ref{minimax}).
\par
A triplet $(\xi,\zeta,\Sigma)$ satisfying (\ref{KKT_minimax1})--(\ref{KKT_minimax3}) is actually obtained via the KKT conditions of the following quadratic semidefinite program (QSDP):
\begin{eqnarray} \label{SQSDPsubproblem0}
\begin{array}{lll}
\displaystyle \mini_{(\xi, \zeta, \Sigma) \in {\cal V}} & \displaystyle \langle \nabla f(x), \xi \rangle + \frac{1}{2} \langle H \xi , \xi \rangle + \frac{\sigma}{2} \Vert \zeta \Vert^{2} + \frac{\sigma}{2} \Vert \Sigma \Vert_{{\rm F}}^{2}
\\
{\rm subject ~ to} & \displaystyle g(x) + \nabla g(x)^{\top} \xi + \sigma (\zeta - y) = 0, 
\\
& \displaystyle X(x) + {\cal A}(x) \xi + \sigma (\Sigma - Z) \succeq O,
\end{array}
\end{eqnarray}
which is further transformed to the following problem by eliminating the variable $\zeta$ via the relation $\zeta = y - \frac{1}{\sigma} \{ g(x) + \nabla g(x)^{\top} \xi \}$:
\begin{eqnarray} \label{SQSDPsubproblem}
\begin{array}{lll}
\displaystyle \mini_{(\xi, \Sigma) \in {\cal W}} & \displaystyle \langle \nabla f(x) - \nabla g(x) s, \xi \rangle + \frac{1}{2} \langle M \xi , \xi \rangle + \frac{\sigma}{2} \Vert \Sigma \Vert_{{\rm F}}^{2}
\\
{\rm subject ~ to} & \displaystyle {\cal A}(x) \xi + \sigma (\Sigma - T) \succeq O,
\end{array}
\end{eqnarray}
where
\begin{eqnarray*}
{\cal W} := {\bf R}^{n} \times {\bf S}^{d}, ~ M := H + \frac{1}{\sigma} \nabla g(x) \nabla g(x)^{\top}, ~ s := y - \frac{1}{\sigma} g(x), ~ T := Z - \frac{1}{\sigma}X(x).
\end{eqnarray*}
We employ problem~(\ref{SQSDPsubproblem}) as a subproblem of the proposed method to generate a search direction, although (\ref{SQSDPsubproblem0}) is equivalent to (\ref{SQSDPsubproblem}) in the sense that both global optimal values are equal if exist. This is because (\ref{SQSDPsubproblem}) is expected to be easier to solve than (\ref{SQSDPsubproblem0}) since the number of variables in (\ref{SQSDPsubproblem}) is less than that of (\ref{SQSDPsubproblem0}).
\par
Note that (\ref{SQSDPsubproblem}) is always strictly feasible, that is, Slater's CQ is satisfied, which can be ensured by substituting $(\xi, \Sigma)=(0, I+T)$. Furthermore, if the objective function of (\ref{SQSDPsubproblem}) is strongly convex, that is, the matrix $M$ is positive definite, problem~(\ref{SQSDPsubproblem}) has a unique optimum. This fact is formally stated in the following proposition.
\begin{proposition} \label{subproblem_pro}
Suppose that $M = H + \frac{1}{\sigma} \nabla g(x) \nabla g(x)^{\top} \succ O$. Then, problem~\eqref{SQSDPsubproblem} has a unique optimal solution.
\end{proposition}
Henceforth, we assume that $M \succ O$ and consider solving QSDP~(\ref{SQSDPsubproblem}) at every iteration. We often call this problem a {\it stabilized QSDP subproblem}. It is derived from the existing subproblem proposed in \cite{Wr98} which was referred to as the stabilized subproblem therein because of its calming effect on a sequence of Lagrange multiplier estimates.
\par
The following proposition connects the merit function $F$ to QSDP~(\ref{SQSDPsubproblem}). Its proof is deferred to Appendix~\ref{appendix_proof}.
\begin{proposition}\label{prop_nablaF}
Assume that $M \succ O$. Then, a unique solution of QSDP~\eqref{SQSDPsubproblem}, say $(\xi^{\ast}, \Sigma^{\ast}) \in {\cal W}$, satisfies $\langle \nabla F(x; \sigma, y, Z), \xi^{\ast} \rangle \leq - \langle M \xi^{\ast}, \xi^{\ast} \rangle - \sigma \Vert \Sigma^{\ast} - [T]_{+} \Vert_{{\rm F}}^{2}$. Moreover, $\nabla F(x; \sigma, y, Z) =0$ if and only if $(\xi^{\ast}, \Sigma^{\ast}) = (0, [T]_{+})$.
\end{proposition}

\begin{remark}
Proposition~{\rm \ref{prop_nablaF}} is a novelty of this paper and is not seen in the existing researches related to the stabilized SQP-type methods. This result motivates us to solve QSDP~\eqref{SQSDPsubproblem} in Mini-$F$-Phase.
\end{remark}

\subsection{An inexact solution of the stabilized QSDP subproblem} \label{section_subpro}
Although it is often impractical or impossible to solve a QSDP subproblem exactly, the existing SQSDP methods such as \cite{CoRa04,FaNoAp02,FrJaVo07,GoRa10,ZhCh16,ZhZh14} establish global convergence by using exact optima of QSDP subproblems, as far as we have investigated. In this section, we study when we may truncate a process of solving QSDP~(\ref{SQSDPsubproblem}).
\par
Let $(\xi^{\ast}, \Sigma^{\ast}) \in {\cal W}$ be the unique optimum of (\ref{SQSDPsubproblem}). Since (\ref{SQSDPsubproblem}) satisfies Slater's CQ as mentioned before, the KKT conditions hold at $(\xi^{\ast}, \Sigma^{\ast})$, namely, there exists some Lagrange multiplier matrix $\Lambda^{\ast}$ such that
\begin{eqnarray}
M \xi^{\ast} + \nabla f(x) - \nabla g(x) s - {\cal A}^{\ast}(x) \Lambda^{\ast} &=& 0, \label{subkkt-1}
\\
\sigma ( \Sigma^{\ast} - \Lambda^{\ast} ) &=& O,\label{subkkt-2}
\\
{\cal A}(x) \xi^{\ast} + \sigma (\Sigma^{\ast} - T) &\succeq& O, \label{subkkt-3}
\\
\Lambda^{\ast} &\succeq& O, \label{subkkt-4}
\\
\langle {\cal A} (x) \xi^{\ast} + \sigma (\Sigma^{\ast} - T), \Lambda^{\ast} \rangle &=& 0. \label{subkkt-5}
\end{eqnarray}
\par
Applying a suitable convergent algorithm such as a primal-dual interior point method to QSDP~(\ref{SQSDPsubproblem}) generates a sequence $\{ (\xi_{j}, \Sigma_{j}, \Lambda_{j}) \}$ which converges to the above $(\xi^{\ast},\Sigma^{\ast},\Lambda^{\ast})$. In particular, this sequence satisfies the following conditions:
\begin{eqnarray}
\begin{array}{ccc} \label{property1-3}
& \hspace{5mm} \displaystyle \lim_{j \to \infty} \xi_{j} = \xi^{\ast}, ~ \lim_{j \to \infty} \Sigma_{j} = \Sigma^{\ast} \succeq O, ~ \lim_{j \to \infty} \Lambda_{j} = \Lambda^{\ast} \succeq O, ~ \lim_{j \to \infty} \eta_{j} = \eta^{\ast} = 0, &
\\
& \hspace{5mm} \displaystyle \lim_{j \to \infty} \Theta_{j} = \Theta^{\ast} = O, ~ \lim_{j \to \infty} \Omega_{j} = \Omega^{\ast} \succeq O, ~  \lim_{j \to \infty} \langle \Omega_{j}, \Lambda_{j} \rangle = \langle \Omega^{\ast}, \Lambda^{\ast} \rangle = 0,&
\end{array}
\end{eqnarray}
where $\{( \eta_{j}, \Theta_{j}, \Omega_{j} )\}$ is a convergent sequence associated with violation error for (\ref{subkkt-1})--(\ref{subkkt-5}) and is defined by
\begin{eqnarray}
\eta_{j} &:=& M \xi_{j} + \nabla f(x) - \nabla g(x) s - {\cal A}^{\ast}(x) \Lambda_{j}, \label{KKT1}
\\
\Theta_{j} &:=& \sigma ( \Sigma_{j} - \Lambda_{j} ), \label{KKT2}
\\
\Omega_{j} &:=& {\cal A}(x) \xi_{j} + \sigma ( \Sigma_{j} - T ) \label{KKT3}
\end{eqnarray}
for each $j \in {\bf N} \cup \{0\}$. As indicated by the convergence analysis in the next section, sufficient decrease of the merit function $F$ with an appropriate $\sigma$ leads to a point satisfying optimality conditions for NSDP~\eqref{NSDP}. Hence, if we find a point $(\xi_{j}, \Sigma_{j}, \Lambda_{j})$ which decreases $F$ sufficiently, we may stop solving QSDP~(\ref{SQSDPsubproblem}) or producing the sequence $\{(\xi_{j}, \Sigma_{j}, \Lambda_{j} )\}$. Specifically, we do it if the following conditions hold for some $j \in {\bf N} \cup \{ 0 \}$:
\begin{eqnarray}
& \langle \nabla F(x; \sigma, y, Z), \xi_{j} \rangle \leq - c_{1} \langle M \xi_{j}, \xi_{j} \rangle - c_{1} \sigma \Vert \Lambda_{j} - [T]_{+} \Vert_{{\rm F}}^{2}, & \label{inexact_ineq1}
\\
& \Vert \eta_{j} \Vert \leq c_{2} | \langle \nabla F(x; \sigma, y, Z), \xi_{j} \rangle |, & \label{inexact_ineq2}
\end{eqnarray}
where $c_{1} \in (0,1)$ and $c_{2} > 0$ are prefixed constants. The above index $j \in {\bf N}$ necessarily exists as is proved in Proposition \ref{inexact_cond}. In the proposed method, we set a search direction as $p := \xi_{j}$, and the next candidates, $\overline{y}$ and $\overline{Z}$, of Lagrange multipliers $y$ and $Z$ as
\begin{eqnarray*}
\overline{y} := \zeta_{j}, \quad \overline{Z} := [\Sigma_{j}]_{+},
\end{eqnarray*}
respectively, where $\zeta_{j} := y - \frac{1}{\sigma} \{ g(x) + \nabla g(x)^{\top} \xi_{j} \}$. One may wonder why we do not determine $y$ and $Z$ immediately in place of $\overline{y}$ and $\overline{Z}$. In the next step, we decide whether or not to set $\overline{y}$ and $\overline{Z}$ to be $y$ and $Z$ in view of some measures related to the optimality conditions.
\par 
Summarizing the above discussion, we can describe the way of solving the QSDP as in Algorithm~\ref{InexatQSDP}.

\begin{remark}
In the field of NLPs, there exists such an existing SQP method which uses an approximate solution of its subproblem. For its details, see \cite{IzSo10}.
\end{remark}

\begin{algorithm}[h]
\caption{Procedure of solving QSDP with a truncation technique} \label{InexatQSDP}
\begin{algorithmic}[1]
\Procedure{{\rm INEXACT-QSDP}}{$\overline{v}_{k+1}, y_{k}, Z_{k}, \phi_{k}, \psi_{k}, \gamma_{k}, \sigma_{k}$}
\State{Let $M := H_{k} + \frac{1}{\sigma_{k}} \nabla g(x_{k}) \nabla g(x_{k})^{\top}, ~ \sigma := \sigma_{k}, ~ x := x_{k}, ~ y := y_{k}, ~ Z := Z_{k}, ~ s := y_{k} - \frac{1}{\sigma_{k}} g(x_{k})$, and $T := Z_{k} - \frac{1}{\sigma_{k}} X(x_{k})$.}
\State{Set $j:=0$}
\Repeat 
\State{Produce $(\xi_{j}, \Sigma_{j}, \Lambda_{j})$ with a suitable algorithm.}
\Comment{The whole sequence $\{ ( \xi_{j}, \Sigma_{j}, \Lambda_{j} ) \}$ is required to satisfies conditions (\ref{property1-3}) related to QSDP~(\ref{SQSDPsubproblem})}
\State{Set $j:=j+1$}
\Until{
$(\xi_{j}, \Sigma_{j}, \Lambda_{j})$ satisfies \eqref{inexact_ineq1} and \eqref{inexact_ineq2}.
} \label{line:stop_tol}
\State{Set $p_{k} := \xi_{j}$, $\overline{y}_{k+1} := y - \frac{1}{\sigma} \{ g(x) + \nabla g(x)^{\top} \xi_{j} \}$, and $\overline{Z}_{k+1} := [\Sigma_{j}]_{+}$.}
\State{{\bf return} $(p_{k},\overline{y}_{k+1},\overline{Z}_{k+1})$} 
\EndProcedure
\end{algorithmic}
\end{algorithm}

\subsection{Line-search along $p_{k}$ and update of $x_{k}$}
From now on, we use a subscript $k \in {\bf N} \cup \{ 0 \}$ to denote a current iterate. After computing the search direction $p_{k}$, we determine a step size along $p_{k}$ so that the function $F(\, \cdot \, ; \sigma_{k}, y_{k}, Z_{k})$ decreases. Notice that (\ref{inexact_ineq1}) with $\xi_{j} = p_{k}$ implies that $\langle \nabla F(x_{k}; \sigma_{k}, y_{k}, Z_{k}), p_{k} \rangle \leq 0$. We adopt a backtracking line-search strategy using $F(\, \cdot \, ; \sigma_{k}, y_{k}, Z_{k})$ to set a step size $\alpha_{k} := \beta^{\ell_{k}}$, where $\beta \in (0,1)$ is a given parameter and $\ell_{k}$ is the smallest nonnegative integer such that
\begin{eqnarray}
F(x_{k} + \beta^{\ell_{k}} p_{k}; \sigma_{k}, y_{k}, Z_{k}) \leq F(x_{k}; \sigma_{k}, y_{k}, Z_{k}) + \tau \beta^{\ell_{k}} \Delta_{k},  \label{eqn:line}
\end{eqnarray}
where $\Delta_{k} := \max \left\{ \langle \nabla F(x_{k}; \sigma_{k}, y_{k}, Z_{k}), p_{k} \rangle, - \omega \Vert p_{k} \Vert^{2} \right\}$, $\omega \in (0,1)$, and $\tau \in (0,1)$. If the value of $| \langle \nabla F(x_{k}; \sigma_{k}, y_{k}, Z_{k}), p_{k} \rangle |$ is large, the term $- \omega \Vert p_{k} \Vert^{2}$ in $\Delta_{k}$ helps us to adopt $\ell_{k}$ at an early stage of the above procedure. After computing $\alpha_{k} = \beta^{\ell_{k}}$, we update $x_{k}$ as $x_{k+1} := x_{k} + \alpha_{k} p_{k}$.

\subsection{Formal statement of the stabilized SQSDP method}
By summarizing the explanation in the previous sections, the proposed stabilized SQSDP method is described as in Algorithm~\ref{algorithm_SQSDP}.
\par
Finally, it may be worth noting the case where $x_{k}$ satisfies $\nabla F(x_{k}; \sigma_{k}, y_{k}, Z_{k}) = 0$ in Step~2 of Algorithm~\ref{algorithm_SQSDP}. In this case, the optimal solution of QSDP~(\ref{SQSDPsubproblem}) is calculated by $(0,[Z_{k} - \frac{1}{\sigma_{k}} X(x_{k})]_{+})$ according to Proposition~\ref{prop_nablaF}. Therefore, $\overline{v}_{k+1} = (x_{k+1}, \overline{y}_{k+1}, \overline{Z}_{k+1})$ computed in Line~\ref{step1} of Algorithm~\ref{algorithm_SQSDP} is identical to the output of $\mbox{INEXACT-QSDP}(H_{k}, x_{k}, y_{k}, Z_{k}, \sigma_{k})$.
\begin{algorithm}[h]
\caption{Inexact and stabilized SQSDP method (main algorithm)} \label{algorithm_SQSDP}
\begin{algorithmic}[1]
\Require 
Choose $v_{0} := (x_{0}, y_{0}, Z_{0})$ such that $Z_{0} \succeq O$. Set
\begin{eqnarray*}
k := 0, ~ \overline{y}_{0} := y_{0}, ~ \overline{Z}_{0} := Z_{0}, ~ \phi_{0} > 0, ~ \psi_{0} > 0, ~ \gamma_{0} > 0, ~ \sigma_{0} > 0.
\end{eqnarray*}
Select $\tau \in (0,1)$, $\omega \in (0,1)$, $c_{1} \in (0,1)$, $c_{2} > 0$, $\beta \in (0,1)$, $\kappa \in (0,1)$, $y_{\max} > 0$, and $z_{\max} > 0$.

\Repeat
\If{$\Vert \nabla F(x_{k}; \sigma_{k}, y_{k}, Z_{k}) \Vert = 0$}
\Comment{Step~1~(Mini-$F$-Phase)}
\State{
\begin{eqnarray*}
\textstyle x_{k+1} := x_{k}, ~ \overline{y}_{k+1} := y_{k} - \frac{1}{\sigma_{k}} g(x_{k+1}),\ \overline{Z}_{k+1} := [Z_{k} - \frac{1}{\sigma_{k}} X(x_{k+1})]_{+}.
\end{eqnarray*}
}
\label{step1}

\Else
\State{Choose $H_{k} \succ O$. Compute $p_{k}$, $\overline{y}_{k+1}$, and $\overline{Z}_{k+1}$ by
\begin{eqnarray*}
\textstyle (p_{k}, \overline{y}_{k+1}, \overline{Z}_{k+1}) := \mbox{INEXACT-QSDP}(H_{k}, x_{k}, y_{k}, Z_{k}, \sigma_{k}).
\end{eqnarray*}
} \label{st:chooH}

\State{Compute the smallest nonnegative integer $\ell_{k}$ such that
\begin{eqnarray*}
& F(x_{k} + \beta^{\ell_{k}} p_{k}; \sigma_{k}, y_{k}, Z_{k}) \leq F(x_{k}; \sigma_{k}, y_{k}, Z_{k}) + \tau \beta^{\ell_{k}} \Delta_{k}, &
\\
& \Delta_{k} := \max \{ \langle \nabla F(x_{k}; \sigma_{k}, y_{k}, Z_{k}), p_{k} \rangle, - \omega \Vert p_{k} \Vert^{2} \}. &
\end{eqnarray*}

\State{Set $x_{k+1} := x_{k} + \beta^{\ell_{k}} p_{k}$.}\label{st:xk}
}
\EndIf
\State{Compute $y_{k+1}$, $Z_{k+1}$, $\phi_{k+1}$, $\psi_{k+1}$, and $\gamma_{k+1}$ by  \Comment{Step~2}
\begin{eqnarray*}
&& (y_{k+1},Z_{k+1},\phi_{k+1},\psi_{k+1},\gamma_{k+1}) 
\\
&& \hspace{20mm} := \mbox{{\rm VOMF-ITERATES}}(\overline{v}_{k+1}, y_{k}, Z_{k}, \phi_{k}, \psi_{k}, \gamma_{k}, \sigma_{k}).
\end{eqnarray*}
}
\State{
Update $\sigma_{k}$ by  \Comment{Step~3}
\begin{eqnarray*} 
&& \sigma_{k+1} := \left\{
\begin{array}{ll}
\min \{ \frac{1}{2} \sigma_{k}, r(v_{k+1})^{\frac{3}{2}} \} & {\rm if ~} \Vert \nabla F(x_{k+1}; \sigma_{k}, y_{k}, Z_{k}) \Vert \leq \gamma_{k},
\\
\sigma_{k} & {\rm otherwise}.
\end{array}
\right.
\end{eqnarray*}}
\State{Set $k := k+1$.}  \Comment{Step~4}
\Until{$v_{k}:=(x_{k}, y_{k}, Z_{k})$ meets a suitable criterion.}
\end{algorithmic}
\end{algorithm}

\section{Convergence Analysis of the stabilized SQSDP method} \label{sec:convergence_analysis}
We show the global convergence of Algorithm~\ref{algorithm_SQSDP}.
\subsection{Well-definedness of Algorithm~\ref{algorithm_SQSDP}} \label{sec:well-definedness}
In what follows, we assume that $M \succ O$ and also that $\{ (\xi_{j}, \Sigma_{j}, \Lambda_{j}, \eta_{j}, \Theta_{j}, \Omega_{j}) \}$ is a sequence produced in INEXACT-QSDP. Recall that this sequence satisfies (\ref{property1-3}). Moreover, we often use
\begin{eqnarray} \label{relation_R}
\begin{array}{c}
\hspace{-30mm} R_{j} := \langle \eta_{j}, \xi_{j} \rangle + \langle \Omega_{j}, \Lambda_{j} \rangle - \langle \Omega_{j}, [T]_{+} \rangle 
\\
\hspace{30mm} + \sigma \langle \Lambda_{j} - [T]_{+}, T - [T]_{+} \rangle - \langle \Lambda_{j} - [T]_{+}, \Theta_{j} \rangle
\end{array}
\end{eqnarray}
for each $j \in {\bf N} \cup \{ 0 \}$.
\par
Our goal in Section~\ref{sec:well-definedness} is to prove that Algorithm~\ref{algorithm_SQSDP} is well-defined in the sense that INEXACT-QSDP is terminated within a finite number of iterations. Specifically, we aim to prove the following proposition:
\begin{pro} \label{inexact_cond}
If $\nabla F(x; \sigma, y, Z) \not = 0$, then there exists $j_{0} \in {\bf N}$ such that \eqref{inexact_ineq1} and \eqref{inexact_ineq2} hold for all $j \geq j_{0}$, meaning that INEXACT-QSDP terminates finitely. 
\end{pro}
In order to prove this proposition, we prepare three lemmas. In the first lemma, we show that $\xi_{j}$ is a descent direction of the merit function $F$ if 
\begin{eqnarray}
R_{j} \leq (1-c_{1}) \langle M \xi_{j}, \xi_{j} \rangle + \sigma( 1 - c_{1} ) \Vert \Lambda_{j} - [T]_{+} \Vert_{{\rm F}}^{2} . \label{suff_cond_descent}
\end{eqnarray}

\begin{lemma} \label{descent_direction}
If there exists $j \in {\bf N}$ satisfying \eqref{suff_cond_descent}, then \eqref{inexact_ineq1} holds.
\end{lemma}

\noindent
{\it Proof.}
Since $\nabla F(x; \sigma, y, Z) = \nabla f(x) - \nabla g(x) s - {\cal A}^{\ast}(x)[T]_{+}$ by (\ref{gradient_F}), we have
\begin{eqnarray}
&& \hspace{-7.5mm} \langle \nabla F(x; \sigma, y, Z), \xi_{j} \rangle \nonumber
\\
&& \hspace{-5mm} = - \langle M \xi_{j}, \xi_{j} \rangle + \langle \eta_{j}, \xi_{j} \rangle + \langle \Lambda_{j} - [T]_{+}, {\cal A}(x) \xi_{j} \rangle  \label{fromKKT1}
\\
&& \hspace{-5mm} = - \langle M \xi_{j}, \xi_{j} \rangle + \langle \eta_{j}, \xi_{j} \rangle + \langle \Lambda_{j} - [T]_{+}, \Omega_{j} \rangle  \label{fromKKT2}
\\
&& \hspace{6mm} + \sigma \langle \Lambda_{j} - [T]_{+}, T - [T]_{+} \rangle + \sigma \langle \Lambda_{j} - [T]_{+}, [T]_{+} - \Sigma_{j} \rangle  \nonumber
\\
&& \hspace{-5mm} = - \langle M \xi_{j}, \xi_{j} \rangle + \langle \eta_{j}, \xi_{j}\rangle + \langle \Omega_{j}, \Lambda_{j} \rangle - \langle \Omega_{j}, [T]_{+} \rangle  \label{fromKKT3}
\\
&& \hspace{6mm} + \sigma \langle \Lambda_{j} - [T]_{+}, T - [T]_{+} \rangle - \sigma \Vert \Lambda_{j} - [T]_{+} \Vert_{{\rm F}}^{2} - \langle \Lambda_{j} - [T]_{+}, \Theta_{j} \rangle  \nonumber
\\
&& \hspace{-5mm} = - \langle M \xi_{j}, \xi_{j} \rangle - \sigma \Vert \Lambda_{j} - [T]_{+} \Vert_{{\rm F}}^{2} + R_{j}, \label{nablaFxi}
\end{eqnarray}
where (\ref{fromKKT1}), (\ref{fromKKT2}), and (\ref{fromKKT3}) are derived from (\ref{KKT1}), (\ref{KKT3}), and (\ref{KKT2}), respectively. The combination of (\ref{suff_cond_descent}) and (\ref{nablaFxi}) yields that $\langle \nabla F(x; \sigma, y, Z), \xi_{j} \rangle \leq - c_{1} \langle M \xi_{j}, \xi_{j} \rangle - c_{1} \sigma \Vert \Lambda_{j} - [T]_{+} \Vert_{{\rm F}}^{2}$. $ \hfill \Box $
\bigskip

In the second lemma, we prove that there exists some index $j \in {\bf N}$ fulfilling condition (\ref{suff_cond_descent}) if the following condition holds:
\begin{eqnarray}
\exists K > 0 \quad {\rm s.t.} \quad K \leq \langle M \xi_{j}, \xi_{j} \rangle + \sigma \Vert \Lambda_{j} - [T]_{+} \Vert_{{\rm F}}^{2} \quad \forall j \in {\bf N}. \label{away_from_K}
\end{eqnarray}

\begin{lemma} \label{suff_descent}
Assume that \eqref{away_from_K} holds. Then, there exists $n_{0} \in {\bf N}$ such that \eqref{inexact_ineq1} holds for all $j \geq n_{0}$.
\end{lemma}

\noindent
{\it Proof.}
It is clear that $\langle \Omega^{\ast}, [T]_{+} \rangle \geq 0$ from $\Omega^{\ast} \succeq O$. Moreover, we see that $\langle \Lambda^{\ast} - [T]_{+}, T - [T]_{+} \rangle \leq 0$ by the facts that $\Lambda^{\ast} \succeq O$ and $[ \, \cdot \, ]_{+}$ is the projection onto the convex set ${\bf S}^{d}_{+}$. We have from these facts, (\ref{property1-3}), and (\ref{relation_R}) that
\begin{eqnarray}
R^{\ast} := \lim_{j \to \infty}R_{j} = - \langle \Omega^{\ast}, [T]_{+} \rangle + \sigma \langle \Lambda^{\ast} - [T]_{+}, T - [T]_{+} \rangle \leq 0. \label{Rstarineq}
\end{eqnarray}
Since (\ref{Rstarineq}) holds, there exists $n_{0} \in {\bf N}$ such that $R_{j} \leq R_{j} - R^{\ast}$ and $| R_{j} - R^{\ast} | \leq (1 - c_{1}) K$ for all $j \geq n_{0}$, that is, $R_{j} \leq (1 - c_{1}) K$ for all $j \geq n_{0}$. Combining this result and (\ref{away_from_K}) yields that (\ref{suff_cond_descent}) holds for all $j \geq n_{0}$. It then follows from Lemma \ref{descent_direction} that (\ref{inexact_ineq1}) holds for all $j \geq n_{0}$.
$ \hfill \Box $
\bigskip

By Lemma~\ref{suff_descent}, we ensure the existence of the index $j$ such that $\xi_{j}$ is the descent direction of $F$ if (\ref{away_from_K}) holds. In the third lemma, we show that $\nabla F(x; \sigma, y, Z) \not = 0$ is a sufficient condition for (\ref{away_from_K}).

\begin{lemma} \label{descent_direction2}
If $\nabla F(x; \sigma, y, Z) \not = 0$, then \eqref{away_from_K} holds.
\end{lemma}

\noindent
{\it Proof.}
We prove the assertion by contradiction. To this end, assume that there exists ${\cal L} \subset {\bf N}$ such that
\begin{eqnarray*}
\langle M \xi_{j}, \xi_{j} \rangle + \sigma \Vert \Lambda_{j} - [T]_{+} \Vert_{{\rm F}}^{2} \to 0 \quad ({\cal L} \ni j \to 0).
\end{eqnarray*}
Considering $M \succ O$ and (\ref{property1-3}) implies that $\xi^{\ast} = 0$ and $\Lambda^{\ast} = [T]_{+}$. It then follows from (\ref{property1-3}) and (\ref{KKT1}) that $0 = \lim_{{\cal L} \ni j \to \infty} \eta_{j} = M \xi^{\ast} + \nabla f(x) - \nabla g(x) s - {\cal A}^{\ast}(x) \Lambda^{\ast} = \nabla f(x) - \nabla g(x) s - {\cal A}^{\ast}(x) [T]_{+}$. Now, recall that $s = y - \frac{1}{\sigma}g(x)$ and $T = Z - \frac{1}{\sigma} X(x)$. As a result, we obtain $\nabla F(x; \sigma, y, Z) = 0$ by (\ref{gradient_F}). However, this contradicts $\nabla F(x; \sigma, y, Z) \not = 0$.
$ \hfill \Box $
\bigskip

By using Lemmas \ref{suff_descent} and \ref{descent_direction2}, Proposition \ref{inexact_cond} is proven as follows.
\bigskip

\noindent
{\bf Proof of Proposition \ref{inexact_cond}.}
By Lemmas \ref{suff_descent} and \ref{descent_direction2}, there exists $n_{0} \in {\bf N}$ such that 
\begin{eqnarray}
\langle \nabla F(x; \sigma, y, Z), \xi_{j} \rangle \leq - c_{1} \langle M \xi_{j}, \xi_{j} \rangle - c_{1} \sigma \Vert \Lambda_{j} - [T]_{+} \Vert_{{\rm F}}^{2} \quad \forall j \geq n_{0}. \label{pro_inexact1}
\end{eqnarray}
Lemma \ref{descent_direction2} and (\ref{pro_inexact1}) derive that there exists $K > 0$ such that
\begin{eqnarray}
&& c_{1}c_{2} K \leq c_{1}c_{2} \langle M \xi_{j}, \xi_{j} \rangle + c_{1}c_{2} \sigma \Vert \Lambda_{j} - [T]_{+} \Vert_{{\rm F}}^{2} \nonumber
\\
&& \phantom{c_{1}c_{2} K} \leq - c_{2} \langle \nabla F(x; \sigma, y, Z), \xi_{j} \rangle \quad \forall j \geq n_{0}. \label{away_from_0}
\end{eqnarray}
Since (\ref{property1-3}) indicates that $\Vert \eta_{j} \Vert \to 0$ as $j \to \infty$, there exists $n_{1} \in {\bf N}$ such that $\Vert \eta_{j} \Vert \leq c_{1} c_{2} K$ for all $j \geq n_{1}$. This fact and (\ref{away_from_0}) yield that
\begin{eqnarray}
\Vert \eta_{j} \Vert \leq c_{2} | \langle \nabla F(x; \sigma, y, Z), \xi_{j} \rangle | \quad \forall j \geq \max \{ n_{0}, n_{1} \}. \label{pro_inexact2}
\end{eqnarray}
It follows from (\ref{pro_inexact1}) and (\ref{pro_inexact2}) that there exists $j_{0} := \max \{ n_{0}, n_{1} \}$ such that (\ref{inexact_ineq1}) and (\ref{inexact_ineq2}) hold for all $j \geq j_{0}$. The proof is complete. $ \hfill \Box $


\subsection{Global convergence of Algorithm~\ref{algorithm_SQSDP}} \label{sec:global_convergence}
In what follows, we prove the global convergence property of Algorithm \ref{algorithm_SQSDP}. Notice that Algorithm~\ref{algorithm_SQSDP} is within a scope of the global convergence analysis of Algorithm~\ref{algorithm_SQSDP0} that we have established in Section~\ref{sec:prototype_global_convergence}, which supposes Assumptions~\ref{assumption_global0} and \ref{assumption_IJK}. The aim in this section is to prove that Assumption~\ref{assumption_IJK} actually holds as for Algorithm~\ref{algorithm_SQSDP} under the presence of Assumption~\ref{assumption_global0} together with the following one:
\begin{Ass}\label{assumption_global4}
There exist positive constants $\nu_{1}$ and $\nu_{2}$ such that 
\begin{eqnarray*}
\begin{array}{c}
\nu_{1} \leq \lambda_{\min}( H_{k} + \frac{1}{\sigma_{k}} \nabla g(x_{k}) \nabla g(x_{k})^{\top} ), \quad \lambda_{\max} (H_{k}) \leq \nu_{2},
\end{array}
\end{eqnarray*}
for all $k \in {\bf N} \cup \{ 0 \}$.  
\end{Ass}
This assumption is controllable in the sense that it certainly holds if we choose $H_{k}$ to be a bounded positive definite matrix for each $k \in {\bf N} \cup \{ 0 \}$. Indeed, for $a$ and $b$ with $0 < a < b$, we can set $H_{k}$ so that $a I \preceq H_{k} \preceq b I$ for each $k$, e.g., $H_{k} = \frac{a+b}{2} I$. Then, the above assumption holds with $(\nu_1,\nu_2)=(a,b)$.
\par
As a blanket assumption, we suppose that Algorithm~\ref{algorithm_SQSDP} generates an infinite set of iterations.
\par
For simplicity, we denote 
\begin{eqnarray}
M_{k} &:=& H_{k} + \frac{1}{\sigma_{k}} \nabla g(x_{k}) \nabla g(x_{k})^{\top}, \label{definition_M}
\\
s_{k} &:=& y_{k} - \frac{1}{\sigma_{k}} g(x_{k}), \label{definition_s}
\\
T_{k} &:=& Z_{k} - \frac{1}{\sigma_{k}} X(x_{k}). \label{definition_T}
\end{eqnarray}
We also use the following notation:
\begin{eqnarray}
\eta_{k} &:=& M_{k} p_{k} + \nabla  f(x_{k}) - \nabla g(x_{k}) s_{k} - {\cal A}^{\ast}(x_{k}) \Lambda_{k}, \label{q1def}
\\
\Theta_{k} &:=& \sigma_{k} ( \Sigma_{k} - \Lambda_{k} ), \nonumber
\\
\Omega_{k} &:=& {\cal A}(x_{k})p_{k} + \sigma_{k} ( \Sigma_{k} - T_{k} ), \nonumber
\end{eqnarray}
where $\Sigma_{k}$ and $\Lambda_{k}$ are final iteration points of the finite sequences $\{ \Sigma_{j} \}$ and $\{ \Lambda_{j} \}$, respectively, generated in INEXACT-QSDP. From Line~\ref{line:stop_tol} in INEXACT-QSDP, it is clear that
\begin{eqnarray}
& \hspace{-5mm} \langle \nabla F(x_{k}; \sigma_{k}, y_{k}, Z_{k}), p_{k} \rangle \leq -c_{1} \langle M_{k} p_{k}, p_{k} \rangle - c_{1} \sigma_{k} \Vert \Lambda_{k} - [T_{k}]_{+} \Vert_{{\rm F}}^{2} \leq 0, & \label{step32-1}
\\
& \hspace{-5mm} \Vert \eta_{k} \Vert \leq c_{2} | \langle \nabla F(x_{k}; \sigma_{k}, y_{k}, Z_{k}), p_{k} \rangle |. & \label{step32-2}
\end{eqnarray}

From now, we will show that there never occurs a situation that ${\rm card}({\cal I}) < \infty$, ${\rm card}({\cal J}) < \infty$, and ${\rm card}({\cal K}) = \infty$, where recall that the sets ${\cal I}$, ${\cal J}$, and ${\cal K}$ are related to VOMF-ITERATES presented as Algorithm~\ref{ProcPara} and are defined by \eqref{def:sets_IJK}. To this end, we prepare the following lemma. Its proof is given in Appendix~\ref{app:lem3}.

\begin{lemma} \label{Kinfinity}
Suppose that Assumptions~{\rm \ref{assumption_global0}} and {\rm \ref{assumption_global4}} hold. If there occurs a situation that ${\rm card}({\cal I}) < \infty$, ${\rm card}({\cal J}) < \infty$, and ${\rm card}({\cal K}) = \infty$, then
\begin{description}
\item[{\rm (i)}] there exist $k_{0} \in {\bf N}$, $\widehat{y} \in {\bf R}^{m}$, $\widehat{Z} \in {\bf S}^{d}$, $\widehat{\gamma} \in {\bf R}$, and $\widehat{\sigma} \in {\bf R}$ such that  $k \in {\cal K}$, $y_{k} = \widehat{y}$, $Z_{k} = \widehat{Z}$, $\gamma_{k} = \widehat{\gamma}$, and $\sigma_{k} = \widehat{\sigma}$ for all $k \geq k_{0}$;
\item[{\rm (ii)}] $\{ p_{k} \}_{k \geq k_{0}}$ is bounded;
\item[{\rm (iii)}] $\liminf_{k \to \infty} |\Delta_{k}| > 0$.
\end{description}
\end{lemma}

Using Lemma~\ref{Kinfinity}, we show that Algorithm~\ref{algorithm_SQSDP} does not generate an infinite set of iterations satisfying ${\rm card}({\cal I}) < \infty, ~ {\rm card}({\cal J}) < \infty$, and ${\rm card}({\cal K}) = \infty$.

\begin{theorem} \label{Kfinite}
Suppose that Assumptions~{\rm \ref{assumption_global0}} and {\rm \ref{assumption_global4}} hold. Then, there never occurs a situation that ${\rm card}({\cal I}) < \infty$, ${\rm card}({\cal J}) < \infty$, and ${\rm card}({\cal K}) = \infty$, meaning that Assumption~{\rm \ref{assumption_IJK}} holds.
\end{theorem}

\noindent
{\it Proof.}
We prove this theorem by contradiction. Suppose that Algorithm \ref{algorithm_SQSDP} generates an infinite set of iterations satisfying ${\rm card}({\cal I}) < \infty, ~ {\rm card}({\cal J}) < \infty$, and ${\rm card}({\cal K}) = \infty$. By Lemma \ref{Kinfinity} (i), there exists $k_{0} \in {\bf N}$ such that $k \in {\cal K}$, $y_{k} = \widehat{y}$, $Z_{k} = \widehat{Z}$, $\gamma_{k} = \widehat{\gamma}$, and $\sigma_{k} = \widehat{\sigma}$ for all $k \geq k_{0}$. From now on, we assume that $k \geq k_{0}$. Notice that the if-statement of Step~1 is false because if there exists an iteration $k_{1}$ such that it is true, then $k_{1}$ is included in ${\cal J}$ at least, that is, it contradicts $k_{1} \in {\cal K}$. In what follows, we derive a contradiction with Lemma \ref{Kinfinity} (iii). Note that $\Delta_{k} \leq 0$ due to (\ref{step32-1}). We have from \eqref{eqn:line} that $F(x_{k+1}; \widehat{\sigma},\widehat{y}, \widehat{Z}) \leq F(x_{k}; \widehat{\sigma}, \widehat{y}, \widehat{Z}) + \tau \beta^{\ell_{k}} \Delta_{k}$, that is,
\begin{eqnarray}
0 \leq - \tau \beta^{\ell_{k}} \Delta_{k} \leq F(x_{k}; \widehat{\sigma}, \widehat{y}, \widehat{Z}) - F(x_{k+1}; \widehat{\sigma}, \widehat{y}, \widehat{Z}). \label{meritF_decrease}
\end{eqnarray}
By Assumption~\ref{assumption_global0}, it is clear that $\{ F(x_{k}; \widehat{\sigma}, \widehat{y}, \widehat{Z}) \}_{k \geq k_{0}}$ is bounded below. Furthermore, it is a non-increasing sequence. These facts and (\ref{meritF_decrease}) yield $\beta^{\ell_{k}} \Delta_{k} \to 0 ~ (k \to \infty)$. Then, there are the two cases: (a) $\lim \inf_{k \to \infty} \beta^{\ell_{k}} > 0$; (b) $\lim \inf_{k \to \infty} \beta^{\ell_{k}} = 0$.
\\
{\bf Case (a):} We readily have
\begin{eqnarray}
\Delta_{k} \to 0 \quad (k \to \infty). \label{case_a1}
\end{eqnarray}
{\bf Case (b):} In this case, there exists ${\cal L}_{1} \subset {\bf N}$ such that $\ell_{k} \to \infty$ as ${\cal L}_{1} \ni k \to \infty$. Moreover, Assumption \ref{assumption_global0} (A2) implies that there exist ${\cal L}_{2} \subset {\bf N}$ and $x^{\ast} \in \Gamma$ such that $x_{k} \to x^{\ast}$ as ${\cal L}_{2} \ni k \to \infty$. Let ${\cal L} := {\cal L}_{1} \cap {\cal L}_{2}$, $k \in {\cal L}$, and $\delta_{k} := \beta^{\ell_{k}-1} (>0)$. Without loss of generality, we assume $\ell_{k} \geq 1$ for all $k \in {\cal L}$ because $\ell_{k} \to \infty$ as ${\cal L} \ni k \to \infty$. Recall that $\ell_{k}$ is the smallest nonnegative integer such that $F(x_{k} + \beta^{\ell_{k}} p_{k}; \widehat{\sigma}, \widehat{y}, \widehat{Z}) \leq F(x_{k}; \widehat{\sigma}, \widehat{y}, \widehat{Z}) + \tau \beta^{\ell_{k}} \Delta_{k}$. Since $\ell_{k}-1$ does not satisfy the inequality, we get $F(x_{k} + \delta_{k} p_{k}; \widehat{\sigma}, \widehat{y}, \widehat{Z}) > F(x_{k}; \widehat{\sigma}, \widehat{y}, \widehat{Z}) + \tau \delta_{k} \Delta_{k}$. From $\langle \nabla F(x_{k}; \widehat{\sigma}, \widehat{y}, \widehat{Z}), p_{k} \rangle \leq \max \{ \langle \nabla F(x_{k}; \widehat{\sigma}, \widehat{y}, \widehat{Z}), p_{k} \rangle, - \omega \Vert p_{k} \Vert^{2} \} = \Delta_{k} \leq 0$, we have
\begin{eqnarray*}
(\tau - 1) \Delta_{k} < \frac{F(x_{k} + \delta_{k} p_{k}; \widehat{\sigma}, \widehat{y}, \widehat{Z}) - F(x_{k}; \widehat{\sigma}, \widehat{y}, \widehat{Z})}{\delta_{k}} - \langle \nabla F(x_{k}; \widehat{\sigma}, \widehat{y}, \widehat{Z}), p_{k} \rangle.
\end{eqnarray*}
Applying the mean value theorem derives that there exists $\theta_{k} \in (0,1)$ satisfying
\begin{eqnarray}
0 \leq (\tau - 1) \Delta_{k} < \langle \nabla F(x_{k} + \theta_{k} \delta_{k} p_{k}; \widehat{\sigma}, \widehat{y}, \widehat{Z}) - \nabla F(x_{k}; \widehat{\sigma}, \widehat{y}, \widehat{Z}), p_{k} \rangle, \label{Delta_to_0}
\end{eqnarray}
where note that $\tau \in (0,1)$. Notice that $\{ p_{k} \}_{k \geq k_{0}}$ is bounded from Lemma \ref{Kinfinity} (ii). Notice also that $\delta_{k} = \beta^{\ell_{k}-1} \to 0$ as ${\cal L} \ni k \to \infty$. As a result, we see that $x^{\ast} = \lim_{{\cal L} \ni k \to \infty} x_{k} = \lim_{{\cal L} \ni k \to \infty} (x_{k} + \theta_{k} \delta_{k} p_{k})$. Thus, the continuity of $\nabla F$ and (\ref{Delta_to_0}) yield that
\begin{eqnarray}
\Delta_{k} \to 0 \quad ({\cal L} \ni k \to \infty). \label{case_a2}
\end{eqnarray}
We have from (\ref{case_a1}) and (\ref{case_a2}) that $\liminf_{k \to \infty} | \Delta_{k} | = 0$. Therefore, this contradicts Lemma \ref{Kinfinity} (iii).
$ \hfill \Box $
\bigskip

Combining Theorems~\ref{global_convergence}, \ref{global_convergence2}, and \ref{Kfinite} readily yields the following theorem.
\begin{theorem} \label{global_convergence_last}
Suppose that Assumptions~{\rm \ref{assumption_global0}} and {\rm \ref{assumption_global4}} hold. Any accumulation point of $\{ x_{k} \}$, say $x^{\ast}$, satisfies at least one of the following statements:
\begin{description}
\item[{\rm (i)}] $x^{\ast}$ is a TAKKT point of \eqref{NSDP};
\item[{\rm (ii)}] $x^{\ast}$ is an AKKT point of \eqref{NSDP};
\item[{\rm (iii)}] $x^{\ast}$ is an infeasible point of \eqref{NSDP}, but a stationary point of the following optimization problem:
\begin{eqnarray*}
\begin{array}{ll}
\displaystyle \mini_{x \in {\bf R}^{n}} & h(x) := \displaystyle \frac{1}{2} \Vert g(x) \Vert^{2} + \frac{1}{2} \Vert [-X(x)]_{+} \Vert_{{\rm F}}^{2},
\end{array}
\end{eqnarray*}
that is to say, $\nabla h(x^{\ast}) = 0$.
\end{description}
Moreover, if any accumulation point of $\{ x_{k} \}$ satisfies \eqref{MFCQ1} and \eqref{MFCQ2} in the definition of the MFCQ, then $x^{\ast}$ is nothing but a KKT point.
\end{theorem}

\begin{remark}
Section~\ref{sec:global_convergence} analyzes the convergence property regarding an AKKT or a TAKKT point. These results cannot be obtained from a simple generalization of existing researches related to the stabilized SQP-type methods such as \cite{GiRo13,Wr98} because they consider only the convergence to a KKT point. We therefore emphasize that the analysis of the current paper is not just a word-for-word extension of that in the existing stabilized SQP-type methods although Algorithm~\ref{algorithm_SQSDP} indeed utilizes ideas and techniques analogous to \cite{GiRo13,Wr98}.
\end{remark}

\section{Numerical experiments} \label{sec_numerical}
We conduct some numerical experiments to examine the efficiency of Algorithm~\ref{algorithm_SQSDP}. The experiments consist of two parts. In the first part, we solve two degenerate problems such that Slater's CQ does not hold. In the second one, we solve non-degenerate problems satisfying Slater's CQ. For the sake of comparison, we also implemented the augmented Lagrangian (AL) method proposed in \cite{AnHaVi18}. The program was implemented with MATLAB R2020b and ran on a machine with an Intel Core i9-9900K 3.60GHz CPU and 128GB RAM. 
\par
In what follows, we describe the concrete setting of Algorithm~\ref{algorithm_SQSDP} and the AL method. We first explain the setting of Algorithm~\ref{algorithm_SQSDP}. We adopted the following as stopping conditions:
\begin{eqnarray*}
r(v_{k}) := r_{V}(x_{k}) + r_{O}(v_{k}) \leq 10^{-6}, ~ \gamma_{k} \leq 10^{-6}, ~ {\rm or} ~ k = 100,
\end{eqnarray*}
where we notice that the functions $r_{V}$ and $r_{O}$ are defined by (\ref{def:rVrO}). In Step~2, we relaxed the condition $\Vert \nabla F(x_{k}; \sigma_{k}, y_{k}, Z_{k}) \Vert = 0$ to
\begin{eqnarray*}
\Vert \nabla F(x_{k}; \sigma_{k}, y_{k}, Z_{k}) \Vert \leq 10^{-6}.
\end{eqnarray*}
We utilized SDPT3 version 4.0 \cite{ToToTu99,TuToTo03} for solving subproblem~(\ref{SQSDPsubproblem}), where we set its stopping criterion to \verb|gaptol|=$10^{-10}$. To ensure the positive definiteness of $M_{k}$ defined by \eqref{definition_M}, we modified it by
\begin{eqnarray*}
M_{k} \leftarrow \left\{
\begin{array}{ll}
M_{k} & \mbox{if the Cholesky decomposition of $M_{k}$ is successful,}
\\
M_{k} + \mu_{k}I & \mbox{otherwise,}
\end{array}
\right.
\end{eqnarray*}
where $\mu_{k} := |\lambda_{\min}(M_{k})| + 10^{-5}$.
The parameters in Algorithm~\ref{algorithm_SQSDP} were set as follows:
\begin{eqnarray*}
& \tau := 10^{-4}, ~ \omega := 10^{-4}, ~ \beta := 0.5, ~ \kappa := 10^{-5}, ~ y_{\max} := 10^{6}, ~ z_{\max} := 10^{6}, &
\\
& \phi_{0} := 10^{3}, ~ \psi_{0} := 10^{3}, ~ \gamma_{0} := 10^{-1}, ~ \sigma_{0} := 10^{-1}. &
\end{eqnarray*}
Moreover, its initial point was selected as $(x_{0}, y_{0}, Z_{0}) := (0,0,O)$.
\par
Next, we explain the setting of the AL method. The overall AL method is described in Appendix~\ref{appendix_numerical}. We used the MATLAB unconstrained optimizer \verb|fminunc| to find an approximate minimizer $x_{k+1}$ of the unconstrained optimization subproblem in Step~2 which satisfies
\begin{eqnarray*}
\Vert \nabla_{x} L(x_{k+1}, \overline{y}_{k} - \rho_{k} g(x_{k+1}), [ \overline{Z}_{k} - \rho_{k} X(x_{k+1}) ]_{+}) \Vert \leq 10^{-10}.
\end{eqnarray*}
The parameters were set as follows:
\begin{eqnarray*}
\varepsilon := 10^{-6}, ~ \tau := 0.5, ~ \gamma := 2, ~ k_{\max} := 100, ~ y_{\max} := 10^{6}, ~ z_{\max} := 10^{6}, ~ \rho_{0} := 10. 
\end{eqnarray*}
Its initial point $(x_{0}, y_{0}, Z_{0})$ was selected as the same one used in Algorithm~\ref{algorithm_SQSDP}.

\subsection{Solving degenerate problems}
This section deals with two degenerate problems described in \cite[Section 3]{GrRe02}. The first one is of the following form:
\begin{eqnarray} \label{testP1}
\begin{array}{ll}
\displaystyle \mini_{X \in {\bf S}^{N}} & \displaystyle \langle C, X \rangle
\\
{\rm subject ~ to} & \displaystyle  [X]_{jj} = 1, ~ j = 1, \ldots, N,
\\
& \displaystyle  \langle J, X \rangle = 0, ~ X \succeq O,
\end{array}
\end{eqnarray}
where $C \in {\bf S}^{N}$ is a constant matrix whose elements are chosen randomly from $[-1,1]$, and $J := e e^{\top}$. Due to $e^{\top} X e = \langle J, X \rangle = 0$, problem~(\ref{testP1}) has no strictly feasible point. Thus, Slater's CQ fails for this problem.
\par
The second problem is as follows:
\begin{eqnarray} \label{testP2}
\begin{array}{ll}
\displaystyle \mini_{X \in {\bf S}^{N}} & \displaystyle \sum_{j=1}^{N} [\alpha]_{j} \langle A_{j}, X \rangle
\\
{\rm subject ~ to} & \displaystyle \langle A_{j}, X \rangle = [b]_{j}, ~ j = 1, \ldots, M,
\\
& X \succeq O,
\end{array}
\end{eqnarray}
where $M, \, N \in {\bf N}$ are constants with $0 < M \leq N$, and $b:=[0,1,1,\ldots,1]^{\top}\in {\bf R}^{M}$, and $\alpha \in {\bf R}^{N}$ is a constant vector whose $n$-th element is set to zero and the others are chosen randomly from $[0,1]$, and $A_{j} := v_{j} v_{j}^{\top}\in {\bf S}^{N}$ for each $j=1, \ldots, N$, where $v_{1}, \ldots, v_{N} \in {\bf R}^{N}$ are arbitrary orthonormal basis vectors. Problem~(\ref{testP2}) also has no strictly feasible point because $v_{1}^{\top} X v_{1} = \langle A_{1}, X \rangle = [b]_{1} = 0$. Hence, Slater's CQ does not hold.
\par
We generated 10 instances for each of problems~\eqref{testP1} and \eqref{testP2} with $N=5,10,15$, and $20$ in the above manner and applied the proposed SQSDP methods and the AL method to them. Tables~\ref{result1-1}--\ref{result2-2} summarize the obtained results, where ``Averaged $\sharp$ iteration'', ``Averaged time'', ``Averaged $r(v^{\ast})$'', and ``Averaged $\max\{ \Vert y^{\ast} \Vert, \Vert Z^{\ast} \Vert_{{\rm F}} \}$'' represent the averaged results for 10 instances for each $N$, and herein $v^{\ast}$, $y^{\ast}$, and $Z^{\ast}$ express the final iteration point of the sequence $\{ v_{k} \}$, $\{ y_{k} \}$, and $\{ Z_{k} \}$ generated by each algorithm, respectively. $\sharp$ite and cputime(s) stand for the number of iterations and spent computational time in seconds, respectively. In addition, the maximum and minimum of ``$r(v^{\ast})$'' and ``$\max\{ \Vert y^{\ast} \Vert, \Vert Z^{\ast} \Vert_{{\rm F}} \}$'' among the 10 runs are shown in the tables. Tables~\ref{result1-1} and \ref{result1-2} show the results of problem~\eqref{testP1}, while Tables~\ref{result2-1} and \ref{result2-2} do those of problem~\eqref{testP2}.
\par
From the tables, we have the following observations: According to Tables~\ref{result1-1} and \ref{result1-2}, Algorithm~\ref{algorithm_SQSDP} and the AL method stopped by reaching the maximum number of iterations $100$ for all the instances, but the results of the two algorithms are crucially different. Indeed, while the AL method diverged in each $N$, all the final iteration points found by Algorithm~\ref{algorithm_SQSDP} satisfied $r(v^{\ast}) \leq 10^{-3}$, which means that Algorithm~\ref{algorithm_SQSDP} found points satisfying the KKT conditions approximately for all $N$.
\par
Next, according to Tables~\ref{result2-1} and \ref{result2-2} concerning problem~(\ref{testP2}), Algorithm~\ref{algorithm_SQSDP} successfully found KKT points for all the instances, whereas the AL method failed in most cases. Furthermore, from the values of $\max \{ \Vert y^{\ast} \Vert, \Vert Z^{\ast} \Vert_{{\rm F}} \}$ in Tables~\ref{result1-2} and \ref{result2-2}, the sequences of Lagrange multipliers computed by the AL method were likely to diverge. As results unseen from the tables, the V-, O-, M-, and F-iterate are performed for about $33.3\%$, $20.8\%$, $0.3\%$, and $45.6\%$ of all the iterations per a single run of Algorithm~\ref{algorithm_SQSDP} for solving problem~\eqref{testP1}. On the other hand, for problem~\eqref{testP2}, we observed that the V-iterate was employed in almost all the iterations.

\subsection{Solving non-degenerate problems}
Next, we examine the efficiency of Algorithm~\ref{algorithm_SQSDP} for two kinds of non-degenerate problems. Those problems are also solved in \cite{YaYa15,YaYaHa12}. The first one is the Gaussian channel capacity problem:
\begin{eqnarray} \label{testP3}
\begin{array}{ll}
\displaystyle \maxi_{x, t \in {\bf R}^{N}} & \displaystyle \frac{1}{2} \sum_{j=1}^{N} \log (1 + [t]_{j})
\\
{\rm subject ~ to} & \displaystyle \frac{1}{N} \sum_{j=1}^{N} [x]_{j} \leq 1, ~ x \geq 0, ~ t \geq 0,
\\ & \displaystyle \left[
\begin{array}{cc}
1 - a_{j}[t]_{j} & \sqrt{r_{j}}
\\
\sqrt{r_{j}} & a_{j}[x]_{j} + r_{j}
\end{array}
\right] \succeq O, ~ j = 1,\ldots, N,
\end{array}
\end{eqnarray}
where the constants $a_{1}, \ldots a_{N}, \, r_{1}, \ldots, r_{N} \in {\bf R}$ are generated randomly from $[0,1]$.
\par
The second problem is the nearest correlation matrix problem:
\begin{eqnarray} \label{testP4}
\begin{array}{ll}
\displaystyle \mini_{X \in {\bf S}^{N}} & \displaystyle \frac{1}{2} \Vert X - A \Vert_{{\rm F}}^{2}
\\
{\rm subject ~ to} & \displaystyle [X]_{jj} = 1, ~ j = 1,\ldots, N,
\\
& X - \eta I \succeq O,
\end{array}
\end{eqnarray}
where $A \in {\bf S}^{N}$ is a constant matrix, and $\eta \in {\bf R}$ is a positive constant. In the experiments, the elements of $A$ are generated randomly from $[-1,1]$ with $[A]_{jj} = 1$ for $j = 1,\ldots, N$, and we set $\eta := 10^{-3}$. We generate 10 instances of each problem in this manner and apply the algorithms to them. The obtained results are summarized in Tables~\ref{result3-1}--\ref{result4-2}, where each column represents the same as in Tables~\ref{result1-1}--\ref{result2-2}.
\par
Algorithm~\ref{algorithm_SQSDP} used the V-iterates in all the iterations for each problem-instance. From Tables~\ref{result3-1}--\ref{result4-2}, we see that both the algorithms successfully found KKT points for all the problem instances. However, there exist substantial difference between the computational times the algorithms spent. Indeed, from Averaged cputimes(s) of Tables~\ref{result3-1} and \ref{result3-2} on problem~\eqref{testP3}, we see that the AL method was likely to find the KKT points faster than Algorithm~\ref{algorithm_SQSDP}. In contrast, according to Table~\ref{result4-1} and \ref{result4-2} on problem~\eqref{testP4}, Algorithm~\ref{algorithm_SQSDP} tended to be much faster than the AL method.

\begin{table}
\vspace{0mm}
\caption{Performance of Algorithm {\rm \ref{algorithm_SQSDP}} on problem \eqref{testP1}}
\vspace{0mm}
\begin{center} \label{result1-1}
\begin{tabular}{|c|cccc|} 
\hline
$N$ & 5 & 10 & 15 & 20 \\
\hline
$\mbox{Averaged $\sharp$ite}$ & 100.0 & 100.0 & 100.0 & 100.0 \\

$\mbox{Averaged cputime(s)}$ & 17.2 & 35.8 & 118.3 & 554.3 \\

$\mbox{Averaged } r(v^{\ast})$ & 2.4e-03 & 7.0e-03 & 6.6e-03 & 1.5e-02 \\

$\mbox{Max of } r(v^{\ast})$ & 1.5e-02 & 6.7e-02 & 4.2e-02 & 8.2e-02 \\

$\mbox{Min of } r(v^{\ast})$  & 1.5e-04 & 5.4e-05 & 2.3e-05 & 1.6e-04 \\

$\mbox{Averaged } \max\{ \Vert y^{\ast} \Vert, \Vert Z^{\ast} \Vert_{{\rm F}} \}$ & 6.4e+03 & 5.6e+03 & 5.8e+03 & 7.1e+03 \\

$\mbox{Max of } \max\{ \Vert y^{\ast} \Vert, \Vert Z^{\ast} \Vert_{{\rm F}} \}$ & 1.3e+04 & 1.0e+04 & 9.3e+03 & 1.1e+04 \\

$\mbox{Min of } \max\{ \Vert y^{\ast} \Vert, \Vert Z^{\ast} \Vert_{{\rm F}} \}$ & 1.6e+03 & 1.8e+03 & 2.5e+03 & 4.2e+03 \\
\hline
\end{tabular}
\end{center}
\vspace{7mm}
\caption{Performance of the AL method on problem \eqref{testP1}}
\vspace{-3mm}
\begin{center} \label{result1-2}
\begin{tabular}{|c|cccc|} 
\hline
$N$ & 5 & 10 & 15 & 20 \\
\hline
$\mbox{Averaged $\sharp$ite}$ & 100.0 & 100.0 & 100.0 & 100.0 \\

$\mbox{Averaged cputime(s)}$ & 6.0 & 28.7 & 111.7 & 355.6 \\

$\mbox{Averaged } r(v^{\ast})$ & 1.9e+31 & 7.8e+31 & 1.1e+32 & 1.2e+32 \\

$\mbox{Max of } r(v^{\ast})$ & 3.3e+31 & 9.7e+31 & 1.2e+32 & 1.3e+32 \\

$\mbox{Min of } r(v^{\ast})$ & 4.1e+30 & 5.1e+31 & 9.1e+31 & 9.4e+31 \\

$\mbox{Averaged } \max\{ \Vert y^{\ast} \Vert, \Vert Z^{\ast} \Vert_{{\rm F}} \}$ & 8.3e+30 & 1.9e+31 & 2.4e+31 & 2.5e+31 \\

$\mbox{Max of } \max\{ \Vert y^{\ast} \Vert, \Vert Z^{\ast} \Vert_{{\rm F}} \}$ & 1.2e+31 & 2.2e+31 & 2.5e+31 & 2.6e+31 \\

$\mbox{Min of } \max\{ \Vert y^{\ast} \Vert, \Vert Z^{\ast} \Vert_{{\rm F}} \}$ & 3.3e+30 & 1.5e+31 & 2.1e+31 & 2.1e+31 \\
\hline
\end{tabular}
\end{center}
\vspace{7mm}
\caption{Performance of Algorithm {\rm \ref{algorithm_SQSDP}} on problem \eqref{testP2}}
\vspace{-3mm}
\begin{center} \label{result2-1}
\fontsize{6}{0pt}\selectfont
\begin{tabular}{|c|cccccc|} 
\hline
$(N,M)$ & (15,5) & (15,10) & (15,15) & (20,7) & (20,14) & (20,20) \\
\hline
$\mbox{Averaged $\sharp$ite}$ & 3.4 & 4.3 & 3.0 & 4.0 & 4.4 & 3.6 \\

$\mbox{Averaged cputime(s)}$ & 4.1 & 5.1 & 2.9 & 20.7 & 22.8 & 15.8 \\

$\mbox{Averaged } r(v^{\ast})$ & 4.1e-07 & 3.9e-09 & 1.5e-07 & 2.4e-07 & 8.1e-09 & 1.0e-07 \\

$\mbox{Max of } r(v^{\ast})$ & 7.4e-07 & 1.1e-08 & 8.6e-07 & 9.1e-07 & 3.5e-08 & 5.1e-07 \\

$\mbox{Min of } r(v^{\ast})$ & 2.6e-10 & 7.3e-10 & 5.2e-09 & 3.4e-10 & 5.7e-10 & 7.8e-12 \\

$\mbox{Averaged } \max\{ \Vert y^{\ast} \Vert, \Vert Z^{\ast} \Vert_{{\rm F}} \}$ & 2.3e+00 & 2.9e+00 & 2.1e+00 & 2.7e+00 & 2.7e+00 & 2.5e+00 \\

$\mbox{Max of } \max\{ \Vert y^{\ast} \Vert, \Vert Z^{\ast} \Vert_{{\rm F}} \}$ & 3.5e+00 & 3.7e+00 & 2.3e+00 & 3.5e+00 & 3.0e+00 & 3.3e+00 \\

$\mbox{Min of } \max\{ \Vert y^{\ast} \Vert, \Vert Z^{\ast} \Vert_{{\rm F}} \}$ & 1.8e+00 & 2.4e+00 & 1.4e+00 & 2.1e+00 & 2.4e+00 & 2.3e+00 \\
\hline
\end{tabular}

\end{center}
\vspace{7mm}
\caption{Performance of the AL method on problem \eqref{testP2}}
\vspace{-3mm}
\begin{center} \label{result2-2}
\fontsize{6}{0pt}\selectfont
\begin{tabular}{|c|cccccc|} 
\hline
$(N,M)$ & (15,5) & (15,10) & (15,15) & (20,7) & (20,14) & (20,20) \\
\hline
$\mbox{Averaged $\sharp$ite}$ & 100.0 & 100.0 & 100.0 & 100.0 & 100.0 & 100.0 \\

$\mbox{Averaged cputime(s)}$ & 107.3 & 107.3 & 109.6 & 329.7 & 326.7 & 340.0 \\

$\mbox{Averaged } r(v^{\ast})$ & 2.9e+31 & 3.4e+31 & 3.7e+31 & 3.4e+31 & 4.3e+31 & 4.2e+31 \\

$\mbox{Max of } r(v^{\ast})$ & 3.2e+31 & 4.3e+31 & 4.1e+31 & 3.7e+31 & 4.7e+31 & 4.8e+31 \\

$\mbox{Min of } r(v^{\ast})$ & 2.4e+31 & 2.9e+31 & 3.4e+31 & 3.0e+31 & 3.9e+31 & 3.8e+31 \\

$\mbox{Averaged } \max\{ \Vert y^{\ast} \Vert, \Vert Z^{\ast} \Vert_{{\rm F}} \}$ & 1.1e+31 & 1.2e+31 & 1.2e+31 & 1.2e+31 & 1.3e+31 & 1.3e+31 \\

$\mbox{Max of } \max\{ \Vert y^{\ast} \Vert, \Vert Z^{\ast} \Vert_{{\rm F}} \}$ & 1.1e+31 & 1.3e+31 & 1.3e+31 & 1.2e+31 & 1.4e+31 & 1.4e+31 \\

$\mbox{Min of } \max\{ \Vert y^{\ast} \Vert, \Vert Z^{\ast} \Vert_{{\rm F}} \}$ & 9.4e+30 & 1.1e+31 & 1.2e+31 & 1.1e+31 & 1.3e+31 & 1.2e+31 \\
\hline
\end{tabular}

\end{center}
\end{table}

\begin{table}
\vspace{0mm}
\caption{Performance of Algorithm {\rm \ref{algorithm_SQSDP}} on problem \eqref{testP3}}
\vspace{-2mm}
\begin{center} \label{result3-1}
\begin{tabular}{|c|cccc|} 
\hline
$N$ & 5 & 10 & 15 & 20 \\
\hline
$\mbox{Averaged $\sharp$ite}$ & 13.7 & 11.0 & 12.7 & 11.8 \\

$\mbox{Averaged cputime(s)}$ & 6.6 & 65.7 & 656.7 & 3051.7 \\

$\mbox{Averaged } r(v^{\ast})$ & 4.0e-07 & 4.2e-07 & 2.7e-07 & 4.1e-07 \\

$\mbox{Max of } r(v^{\ast})$ & 9.9e-07 & 9.7e-07 & 1.0e-06 & 9.6e-07 \\

$\mbox{Min of } r(v^{\ast})$ & 3.1e-08 & 1.5e-07 & 1.7e-08 & 2.3e-08 \\

$\mbox{Averaged } \max\{ \Vert y^{\ast} \Vert, \Vert Z^{\ast} \Vert_{{\rm F}} \}$ & 1.1e+00 & 1.6e+00 & 1.9e+00 & 2.2e+00 \\

$\mbox{Max of } \max\{ \Vert y^{\ast} \Vert, \Vert Z^{\ast} \Vert_{{\rm F}} \}$ & 1.1e+00 & 1.6e+00 & 1.9e+00 & 2.2e+00 \\

$\mbox{Min of } \max\{ \Vert y^{\ast} \Vert, \Vert Z^{\ast} \Vert_{{\rm F}} \}$ & 1.1e+00 & 1.6e+00 & 1.9e+00 & 2.2e+00 \\
\hline
\end{tabular}
\end{center}
\vspace{7mm}
\caption{Performance of the AL method on problem \eqref{testP3}}
\vspace{-2mm}
\begin{center} \label{result3-2}
\begin{tabular}{|c|cccc|} 
\hline
$N$ & 5 & 10 & 15 & 20 \\
\hline
$\mbox{Averaged $\sharp$ite}$ & 6.0 & 6.0 & 6.5 & 6.3 \\

$\mbox{Averaged cputime(s)}$ & 0.2 & 0.2 & 0.7 & 1.5 \\

$\mbox{Averaged } r(v^{\ast})$ & 5.6e-08 & 1.1e-07 & 3.0e-07 & 3.1e-07 \\

$\mbox{Max of } r(v^{\ast})$ & 5.6e-08 & 1.1e-07 & 9.0e-07 & 1.0e-06 \\

$\mbox{Min of } r(v^{\ast})$ & 5.6e-08 & 1.1e-07 & 1.0e-07 & 1.1e-08 \\

$\mbox{Averaged } \max\{ \Vert y^{\ast} \Vert, \Vert Z^{\ast} \Vert_{{\rm F}} \}$ & 1.1e+00 & 1.6e+00 & 1.9e+00 & 2.2e+00 \\

$\mbox{Max of } \max\{ \Vert y^{\ast} \Vert, \Vert Z^{\ast} \Vert_{{\rm F}} \}$ & 1.1e+00 & 1.6e+00 & 1.9e+00 & 2.2e+00 \\

$\mbox{Min of } \max\{ \Vert y^{\ast} \Vert, \Vert Z^{\ast} \Vert_{{\rm F}} \}$ & 1.1e+00 & 1.6e+00 & 1.9e+00 & 2.2e+00 \\
\hline
\end{tabular}
\end{center}
\vspace{7mm}
\caption{Performance of Algorithm {\rm \ref{algorithm_SQSDP}} on problem \eqref{testP4}}
\vspace{-2mm}
\begin{center} \label{result4-1}
\begin{tabular}{|c|cccc|} 
\hline
$N$ & 5 & 10 & 15 & 20 \\
\hline
$\mbox{Averaged $\sharp$ite}$ & 4.2 & 5.6 & 5.0 & 5.8 \\

$\mbox{Averaged cputime(s)}$ & 0.3 & 1.1 & 2.5 & 10.1 \\

$\mbox{Averaged } r(v^{\ast})$ & 3.8e-07 & 1.6e-07 & 7.0e-08 & 5.2e-08 \\

$\mbox{Max of } r(v^{\ast})$ & 5.5e-07 & 6.8e-07 & 4.1e-07 & 4.0e-07 \\

$\mbox{Min of } r(v^{\ast})$ & 8.6e-08 & 3.5e-09 & 3.2e-09 & 6.5e-12 \\

$\mbox{Averaged } \max\{ \Vert y^{\ast} \Vert, \Vert Z^{\ast} \Vert_{{\rm F}} \}$ & 6.6e-01 & 2.9e+00 & 5.7e+00 & 8.5e+00 \\

$\mbox{Max of } \max\{ \Vert y^{\ast} \Vert, \Vert Z^{\ast} \Vert_{{\rm F}} \}$ & 1.5e+00 & 3.8e+00 & 6.5e+00 & 9.5e+00 \\

$\mbox{Min of } \max\{ \Vert y^{\ast} \Vert, \Vert Z^{\ast} \Vert_{{\rm F}} \}$ & 2.2e-01 & 2.1e+00 & 4.9e+00 & 7.6e+00 \\
\hline
\end{tabular}
\end{center}
\vspace{7mm}
\caption{Performance of the AL method on problem \eqref{testP4}}
\vspace{-2mm}
\begin{center} \label{result4-2}
\begin{tabular}{|c|cccc|} 
\hline
$N$ & 5 & 10 & 15 & 20 \\
\hline
$\mbox{Averaged $\sharp$ite}$ & 9.2 & 12.1 & 13.0 & 14.6 \\

$\mbox{Averaged cputime(s)}$ & 0.3 & 7.2 & 78.0 & 470.5 \\

$\mbox{Averaged } r(v^{\ast})$ & 4.7e-07 & 5.1e-07 & 7.1e-07 & 5.0e-07 \\

$\mbox{Max of } r(v^{\ast})$ & 9.9e-07 & 8.8e-07 & 8.6e-07 & 8.9e-07 \\

$\mbox{Min of } r(v^{\ast})$ & 1.5e-07 & 1.7e-07 & 4.4e-07 & 2.5e-07 \\

$\mbox{Averaged } \max\{ \Vert y^{\ast} \Vert, \Vert Z^{\ast} \Vert_{{\rm F}} \}$ & 6.6e-01 & 2.9e+00 & 5.7e+00 & 8.5e+00 \\

$\mbox{Max of } \max\{ \Vert y^{\ast} \Vert, \Vert Z^{\ast} \Vert_{{\rm F}} \}$ & 1.5e+00 & 3.8e+00 & 6.5e+00 & 9.5e+00 \\

$\mbox{Min of } \max\{ \Vert y^{\ast} \Vert, \Vert Z^{\ast} \Vert_{{\rm F}} \}$ & 2.2e-01 & 2.1e+00 & 4.9e+00 & 7.6e+00 \\
\hline
\end{tabular}
\end{center}
\end{table}

\section{Concluding remarks} \label{sec_conclusion}
In this paper, we have proposed a stabilized SQSDP method (Algorithm~\ref{algorithm_SQSDP}) for NSDP~(\ref{NSDP}). In the algorithm, we approximately solve stabilized QSDP subproblems, which are derived from the Lagrange minimax problem of the NSDP, to generate search directions. Those subproblems are always feasible and have global optima. We integrated the V-, O-, M-, and F-iterates, which originate from the stabilized SQP method \cite{GiRo13} for NLPs, into Algorithm~\ref{algorithm_SQSDP}. We showed that, under some mild assumptions, it globally converges to one of the following points: a stationary point concerning the feasibility, an AKKT point, and a TAKKT point. One remarkable point of that analysis is that we need no CQs or boundedness of the produced dual sequence. In the numerical experiments, we observed that Algorithm~\ref{algorithm_SQSDP} could successfully solve many test problems which the AL method failed to solve. From this fact, Algorithm~\ref{algorithm_SQSDP} is promising in comparison with the AL method \cite{AnHaVi18}.
\par
As a future work, it is an interesting direction to study superlinear convergence of Algorithm~\ref{algorithm_SQSDP}.


\appendix

\makeatletter
\setcounter{equation}{0}
\renewcommand{\theequation}{A.\arabic{equation}}
\makeatother

\makeatletter
\setcounter{equation}{0}
\makeatother

\section{Proof of Lemma~\ref{index_lemma}}\label{app_lemindex}
In the following, we show Lemma~\ref{index_lemma}.
\\
\noindent
{\it Proof}
We show item (i). Note that $\{ \phi_{k} \}$ and $\{ \psi_{k} \}$ are non-increasing, and $\phi_{k+1} = \frac{1}{2} \phi_{k}$ or $\psi_{k+1} = \frac{1}{2} \psi_{k}$ for $k \in {\cal I}$. Since ${\rm card}({\cal I}) = \infty$, it is clear that $\phi_{k} \to 0$ or $\psi_{k} \to 0$ as $k \to \infty$.
\par
We next show item (ii). Considering ${\rm card}({\cal I}) < \infty$ yields that
\begin{eqnarray}
\exists n_{0} \in {\bf N} \quad {\rm s.t} \quad k \in {\cal J} \cup {\cal K} \quad \forall k \geq n_{0}. \label{index_JK}
\end{eqnarray}
Let $\widetilde{z} := \max \{ \lambda_{\max}(Z_{0}), \lambda_{\max}(Z_{1}), \ldots, \lambda_{\max}(Z_{n_{0}}), z_{\max} \}$ and $\widetilde{D} := \{ Z \in {\bf S}^{d} \colon O \preceq Z \preceq \widetilde{z} I \}$. We prove $\{ Z_{k} \} \subset \widetilde{D}$ by mathematical induction. If $k=0$, it is clear that $O \preceq Z_{0} \preceq \lambda_{\max}(Z_{0}) I$, that is, $Z_{0} \in \widetilde{D}$. Let $k$ be a nonnegative integer, and let $Z_{k} \in \widetilde{D}$. Now, we consider two cases: $k \geq n_{0}$, $k < n_{0}$. In the first case, we have from (\ref{index_JK}) that if $k \in {\cal J}$, then $Z_{k+1} = \Pi_{D}([ Z_{k} - \frac{1}{\sigma_{k}} X(x_{k+1}) ]_{+}) \in D \subset \widetilde{D}$; if $k \in {\cal K}$, then $Z_{k+1} = Z_{k} \in \widetilde{D}$. As a result, $Z_{k+1} \in \widetilde{D}$. On the other hand, we consider the case where $k < n_{0}$. Noting $1 \leq k+1 \leq n_{0}$ derives that $O \preceq Z_{k+1} \preceq \max \{ \lambda_{\max}(Z_{1}), \lambda_{\max}(Z_{2}), \ldots, \lambda_{\max}(Z_{n_{0}}) \} I$, which implies $Z_{k+1} \in \widetilde{D}$. Therefore, we obtain that $\{ Z_{k} \} \subset \widetilde{D}$. The fact and the boundedness of $\widetilde{D}$ yield the desired result.
\par
Finally, we prove item (iii). Note that $\{ \gamma_{k} \}$ and $\{ \sigma_{k} \}$ are non-increasing. It is easily seen that $\gamma_{k+1} = \frac{1}{2} \gamma_{k}$ for $k \in {\cal J}$. Moreover, we have from (\ref{update_sigma}) that $\sigma_{k+1} = \min \{ \frac{1}{2} \sigma_{k}, r(v_{k+1})^{\frac{3}{2}} \} \leq \frac{1}{2}\sigma_k$ for $k \in {\cal J}$. These facts and ${\rm card}(\cal{J})=\infty$ imply $\gamma_{k} \to 0$ and $\sigma_k\to 0$ as $k \to \infty$. The proof is complete.
$ \hfill \Box $
\bigskip

\section{Proof of Proposition~\ref{prop_nablaF}} \label{appendix_proof}
In Appendix~\ref{appendix_proof}, we give the proof of Proposition \ref{prop_nablaF}.
\\
\noindent
{\it Proof.}
Note that Proposition~\ref{subproblem_pro} and $M \succ O$ imply that problem (\ref{SQSDPsubproblem}) has a unique optimal solution $(\xi^{\ast}, \Sigma^{\ast})$. First, we prove that $\langle \nabla F(x; \sigma, y, Z), \xi^{\ast} \rangle \leq - \langle M \xi^{\ast}, \xi^{\ast} \rangle - \sigma \Vert \Sigma^{\ast} - [T]_{+} \Vert_{{\rm F}}^{2}$. It follows from (\ref{subkkt-1}), (\ref{subkkt-2}), and (\ref{gradient_F}) that
\begin{eqnarray}
\langle \nabla F(x; \sigma, y, Z), \xi^{\ast} \rangle = - \langle M \xi^{\ast}, \xi^{\ast} \rangle + \langle \Sigma^{\ast}, {\cal A}(x) \xi^{\ast} \rangle - \langle [T]_{+}, {\cal A}(x) \xi^{\ast} \rangle. \label{app_eq1}
\end{eqnarray}
By (\ref{subkkt-2}) and (\ref{subkkt-5}), we easily see that $\langle \Sigma^{\ast}, {\cal A}(x) \xi^{\ast} + \sigma( \Sigma^{\ast} - T ) \rangle = 0$, that is, $\langle \Sigma^{\ast}, {\cal A}(x) \xi^{\ast} \rangle = \sigma \langle \Sigma^{\ast}, T - \Sigma^{\ast} \rangle$. Moreover, we have from (\ref{subkkt-3}) and $[T]_{+} \succeq O$ that $0 \leq \langle [T]_{+}, {\cal A}(x) \xi^{\ast} + \sigma (\Sigma^{\ast} - T) \rangle = \langle [T]_{+}, {\cal A}(x) \xi^{\ast} \rangle + \sigma \langle -[T]_{+}, T - \Sigma^{\ast} \rangle$, and hence
\begin{eqnarray}
\qquad \langle \Sigma^{\ast}, {\cal A}(x) \xi^{\ast} \rangle - \langle [T]_{+}, {\cal A}(x) \xi^{\ast} \rangle
&=& \sigma \langle \Sigma^{\ast}-[T]_{+}, T - \Sigma^{\ast} \rangle \nonumber
\\
&=& \sigma \langle \Sigma^{\ast} - [T]_{+}, T - [T]_{+} \rangle - \sigma \Vert \Sigma^{\ast} - [T]_{+} \Vert_{{\rm F}}^{2}. \label{app_ineq2}
\end{eqnarray}
Now note that $\Sigma^{\ast} \succeq O$ by (\ref{subkkt-2}) and (\ref{subkkt-4}). Since $[ \, \cdot \, ]_{+}$ is the projection onto ${\bf S}^{d}_{+}$, we obtain $\langle \Sigma^{\ast} - [T]_{+}, T - [T]_{+} \rangle \leq 0$. This fact and (\ref{app_ineq2}) yield that $\langle \Sigma^{\ast}, {\cal A}(x) \xi^{\ast} \rangle - \langle [T]_{+}, {\cal A}(x) \xi^{\ast} \rangle \leq -\sigma \Vert \Sigma^{\ast} - [T]_{+} \Vert_{{\rm F}}^{2}$. It then follows from (\ref{app_eq1}) that 
\begin{eqnarray}
\langle \nabla F(x; \sigma, y, Z), \xi^{\ast} \rangle \leq - \langle M \xi^{\ast}, \xi^{\ast} \rangle - \sigma \Vert \Sigma^{\ast} - [T]_{+} \Vert_{{\rm F}}^{2}. \label{app_ineq3}
\end{eqnarray}
\par
Now, we assume that $\nabla F(x; \sigma, y, Z) = 0$. We have from (\ref{app_ineq3}) that $\xi^{\ast} = 0$ and $\Sigma^{\ast} = [T]_{+}$. This means that $(0, [T]_{+})$ is a unique solution of problem (\ref{SQSDPsubproblem}).
\par
Conversely, we suppose that $(0, [T]_{+})$ is a unique solution, namely, $(\xi^{\ast}, \Sigma^{\ast}) = (0, [T]_{+})$. By (\ref{subkkt-2}), we obtain $\Lambda^{\ast} = \Sigma^{\ast} = [T]_{+}$. It then follows from (\ref{subkkt-1}) and (\ref{gradient_F}) that $0 = M \xi^{\ast} + \nabla f(x) - \nabla g(x) s - {\cal A}^{\ast}(x) \Lambda^{\ast} = \nabla f(x) - \nabla g(x) s - {\cal A}^{\ast}(x) [T]_{+} = \nabla F(x; \sigma, y, Z)$. The assertion is proven. $ \hfill \Box $

\section{Proof Lemma~\ref{Kinfinity}}\label{app:lem3}
This appendix gives the proof of Lemma~\ref{Kinfinity}.
\\
\noindent
{\it Proof.}
First, we prove item (i). Since ${\rm card}({\cal I}) < \infty$, ${\rm card}({\cal J}) < \infty$, and ${\rm card}({\cal K}) = \infty$, there exists $k_{0} \in {\bf N}$ such that $k \in {\cal K}$, $y_{k} = \widehat{y}$, $Z_{k} = \widehat{Z}$, and $\gamma_{k} = \widehat{\gamma}$ for all $k \geq k_{0}$, where $\widehat{y} := y_{k_{0}}, ~ \widehat{Z} := Z_{k_{0}}$, and $\widehat{\gamma} := \gamma_{k_{0}}$. Moreover (\ref{update_sigma}) yields that $\sigma_{k} = \widehat{\sigma}$ for all $k \geq k_{0}$, where $\widehat{\sigma} := \sigma_{k_{0}}$.
\par
Next, we show item (ii). From now on, assume that $k \geq k_{0}$, that is, item (i) holds. The combination of (\ref{definition_s}), (\ref{definition_T}), and Assumption~\ref{assumption_global0} (A1) and (A2) implies that $\Vert s_{k} \Vert \leq \Vert \widehat{y} \Vert + \frac{1}{\widehat{\sigma}} \sup_{x \in \Gamma} \Vert g(x) \Vert =: s_{\max} < \infty$ and $\Vert [T_{k}]_{+} \Vert_{{\rm F}} \leq \Vert T_{k} \Vert_{{\rm F}} \leq \Vert \widehat{Z} \Vert_{{\rm F}} + \frac{1}{\widehat{\sigma}} \sup_{x \in \Gamma} \Vert X(x) \Vert_{{\rm F}} =: T_{\max} < \infty$. Assumption~\ref{assumption_global4} and (\ref{step32-1}) derive that $\frac{c_{1}}{\nu_{1}} \Vert p_{k} \Vert^{2} \leq - \langle \nabla F(x_{k}; \widehat{\sigma}, \widehat{y}, \widehat{Z}), p_{k} \rangle \leq \Vert \nabla F(x_{k}; \widehat{\sigma}, \widehat{y}, \widehat{Z}) \Vert \Vert p_{k} \Vert$. These results and (\ref{gradient_F}) yield that
\begin{eqnarray*}
\frac{c_{1}}{\nu_{1}} \Vert p_{k} \Vert 
&\leq&
\left\{ \Vert \nabla f(x_{k}) \Vert + \Vert \nabla g(x_{k}) \Vert_{2} \Vert s_{k} \Vert + \sqrt{\sum_{j=1}^{n} \Vert A_{j}(x_{k}) \Vert_{{\rm F}}^{2} \Vert [T]_{+} \Vert^{2}_{{\rm F}}} \right\}
\\
&\leq& \sup_{x \in \Gamma} \left\{ \Vert \nabla f(x) \Vert + s_{\max} \Vert \nabla g(x) \Vert_{2} + T_{\max} \sqrt{\sum_{j=1}^{n} \Vert A_{j}(x) \Vert_{{\rm F}}^{2}} \right\}.
\end{eqnarray*}
Therefore, the boundedness of $\{ p_{k} \}_{k \geq k_{0}}$ is derived from Assumption~\ref{assumption_global0}.
\par
We prove (iii). To show this, we first verify $\liminf_{k \to \infty} | \langle \nabla F(x_{k}; \sigma_{k}, y_{k}, Z_{k}), p_{k} \rangle | > 0$. Now, we assume to the contrary that $\liminf_{k \to \infty} | \langle \nabla F(x_{k}; \sigma_{k}, y_{k}, Z_{k}), p_{k} \rangle | = 0$, that is, there exists ${\cal L} \subset {\bf N}$ such that 
\begin{eqnarray}
| \langle \nabla F(x_{k}; \sigma_{k}, y_{k}, Z_{k}), p_{k} \rangle | \to 0 \quad ({\cal L} \ni k \to \infty). \label{Ftozero_Lnik}
\end{eqnarray}
Let $k$ be a positive integer satisfying $k \geq k_{0}$. Note that item (i) yields $k \in {\cal K}, ~ y_{k} = \widehat{y}, ~ Z_{k} = \widehat{Z}, ~ \gamma_{k} = \widehat{\gamma}$, and $\sigma_{k} = \widehat{\sigma}$. Now, we can easily see that the if-statement of the M-iterate in Algorithm~\ref{ProcPara} is false for $k \geq k_{0}$, namely,
\begin{eqnarray}
\Vert \nabla F(x_{k}; \widehat{\sigma}, \widehat{y}, \widehat{Z}) \Vert > \widehat{\gamma} > 0 \quad \forall k \geq k_{0}+1. \label{bounded_away_from_zero}
\end{eqnarray}
Meanwhile, it follows from (\ref{step32-1}), (\ref{step32-2}), and Assumption~\ref{assumption_global4} that $0 \leq c_{1} \nu_{1} \Vert p_{k} \Vert^{2} + c_{1} \widehat{\sigma} \Vert \Lambda_{k} - [T_{k}]_{+} \Vert_{{\rm F}}^{2} \leq -\langle \nabla F(x_{k}; \widehat{\sigma}, \widehat{y}, \widehat{Z}), p_{k} \rangle$ and $\Vert \eta_{k} \Vert \leq c_{2} | \langle \nabla F(x_{k}; \widehat{\sigma}, \widehat{y}, \widehat{Z}), p_{k} \rangle |$. By considering $| \langle \nabla F(x_{k}; \sigma_{k}, y_{k}, Z_{k}), p_{k} \rangle | \to 0$ as ${\cal L} \ni k \to \infty$, we obtain $\Vert p_{k} \Vert \to 0, ~ \Vert \Lambda_{k} - [T_{k}]_{+} \Vert_{{\rm F}} \to 0$, and $\Vert \eta_{k} \Vert \to 0$. Moreover, we have from (\ref{gradient_F}), (\ref{definition_M}), (\ref{q1def}), and Assumptions~\ref{assumption_global0} and \ref{assumption_global4} that
\begin{eqnarray*}
&& \hspace{-5mm} \Vert \nabla F(x_{k}; \widehat{\sigma}, \widehat{y}, \widehat{Z}) \Vert \leq \Vert \eta_{k} \Vert + \left\{ \nu_{2} + \frac{1}{\widehat{\sigma}} \sup_{x \in \Gamma} \Vert \nabla g(x) \Vert_{2}^{2} \right\} \Vert p_{k} \Vert
\\
&& \hspace{39.2mm} + \left\{ \sup_{x \in \Gamma} \sqrt{ \sum_{j=1}^{n} \Vert A_{j}(x) \Vert_{{\rm F}}^{2} }\right\} \Vert \Lambda_{k} - [T_{k}]_{+} \Vert_{{\rm F}} \to 0
\end{eqnarray*}
as ${\cal L} \ni k \to \infty$. However, this fact contradicts (\ref{bounded_away_from_zero}). As a result, we conclude that $\liminf_{k \to \infty} | \langle \nabla F(x_{k}; \sigma_{k}, y_{k}, Z_{k}), p_{k} \rangle | > 0$. Secondly, we verify that $\liminf_{k \to \infty} \Vert p_{k} \Vert > 0$ by contradiction. Suppose to the contrary that there exists ${\cal M} \subset {\bf N}$ such that $\Vert p_{k} \Vert \to 0$ as ${\cal M} \ni k \to \infty$. Then, it is clear that $\langle \nabla F(x_{k}; \sigma_{k}, y_{k}, Z_{k}), p_{k} \rangle \to 0$ as ${\cal M} \ni k \to \infty$. Since this result corresponds to \eqref{Ftozero_Lnik}, the remainder of the proof can be given in a way similar to that of $\liminf_{k \to \infty} | \langle \nabla F(x_{k}; \sigma_{k}, y_{k}, Z_{k}), p_{k} \rangle | > 0$. Therefore, we obtain $\liminf_{k \to \infty} | \langle \nabla F(x_{k}; \sigma_{k}, y_{k}, Z_{k}), p_{k} \rangle | > 0$ and $\liminf_{k \to \infty} \Vert p_{k} \Vert > 0$. It follows from $\Delta_{k} = \max \{ \langle \nabla F(x_{k}; \sigma_{k}, y_{k}, Z_{k}), p_{k} \rangle, -\omega \Vert p_{k} \Vert^{2} \}$ that $\liminf_{k \to \infty} | \Delta_{k} | > 0$. The proof is complete. $ \hfill \Box $
\bigskip

\makeatletter
\setcounter{equation}{0}
\makeatother

\section{Supplementary for the numerical experiments}  \label{appendix_numerical}
In Appendix~\ref{appendix_numerical}, we give the existing augmented Lagrangian (AL) method \cite{AnHaVi18}.
\\
\\
{\bf The AL method}
\\
{\bf Step 0:} Select $\varepsilon > 0$, $\tau \in (0,1)$, $\gamma > 1$, $k_{\max} \in {\bf N}$, $y_{\max} > 0$, and $z_{\max} > 0$. 
Choose $(x_{0}, y_{0}, Z_{0}) \in {\cal V}$. Set $\rho_{0} > 0$, $\overline{y}_{0} := y_{0}$, $\overline{Z}_{0} := Z_{0}$, and $k:=0$.
\\
{\bf Step 1:}  If $k = k_{\max}$, $r(v_{k}) \leq \varepsilon$, or the following conditions are satisfied, then stop.
\begin{eqnarray*}
& \Vert \nabla_{x} L(v_{k}) \Vert \leq \varepsilon, ~ \Vert g(x_{k}) \Vert + \Vert [-X(x_{k})]_{+} \Vert_{{\rm F}} \leq \varepsilon, ~ \Vert U_{k} - S_{k} \Vert_{{\rm F}} \leq \varepsilon, &
\\
& \lambda_{j}^{U_{k}}(-X(x_{k})) < -\varepsilon ~ \Longrightarrow ~ |\lambda_{j}^{S_{k}}(Z_{k})| \leq \varepsilon, &
\end{eqnarray*}
where $U_{k}$ and $S_{k}$ are orthogonal matrices such that they diagonalize $X(x_{k})$ and $Z_{k}$, respectively, that is to say, $X(x_{k}) = U_{k} {\rm diag}[ \lambda_{1}^{U_{k}}(X(x_{k})), \ldots, \lambda_{d}^{U_{k}}(X(x_{k})) ] U_{k}^{\top}$ and $Z_{k} = S_{k} {\rm diag}[ \lambda_{1}^{S_{k}}(Z_{k}), \ldots, \lambda_{d}^{S_{k}}(Z_{k}) ] S_{k}^{\top}$.
\\
{\bf Step 2:}  Find an approximate minimizer $x_{k+1}$ of the unconstrained optimization problem
\begin{eqnarray*}
\begin{array}{ll}
\displaystyle \mini_{x \in {\bf R}^{n}} & \displaystyle f(x) + \frac{1}{2\rho_{k}} \Vert \overline{y}_{k} - \rho_{k} g(x) \Vert^{2} + \frac{1}{2\rho_{k}} \Vert [\overline{Z}_{k} - \rho_{k} X(x)]_{+} \Vert_{{\rm F}}^{2}.
\end{array}
\end{eqnarray*}
\\
{\bf Step 3:} Update $\rho_{k}$ by
\begin{eqnarray*}
\rho_{k+1} := \left\{
\begin{array}{ll}
\rho_{k} & {\rm if} ~ u_{k+1} \leq \tau u_{k},
\\
\gamma \rho_{k} & {\rm otherwise},
\end{array}
\right. 
\end{eqnarray*}
where $u_{k} := \max \left\{ \Vert g(x_{k}) \Vert, \Vert V_{k} \Vert_{{\rm F}} \right\}$, $V_{k} := [ \frac{1}{\rho_{k}} \overline{Z}_{k} - X(x_{k}) ]_{+} - \frac{1}{\rho_{k}} \overline{Z}_{k}$.
\\
{\bf Step 4:}  Set
\begin{eqnarray*}
& y_{k+1} := \overline{y}_{k} - \rho_{k} g(x_{k+1}), ~ Z_{k+1} := \left[ \overline{Z}_{k} - \rho_{k} X(x_{k+1}) \right]_{+}, &
\\
& \overline{y}_{k+1} := \Pi_{C}(y_{k+1}), ~ \overline{Z}_{k+1} := \Pi_{D}(Z_{k+1}),
\end{eqnarray*}
where $C = \{ y \in {\bf R}^{m} \colon -y_{\max}e \leq y \leq y_{\max} e \}$ and $D = \{ Z \in {\bf S}^{d} \colon O \preceq Z \preceq z_{\max} I \}$. 
\\
{\bf Step 5:}  Set $k := k+1$, and go back to Step 1.

\end{document}